\documentclass[11pt]{amsart}
\usepackage{color}
\usepackage{epsfig}
\usepackage{verbatim}

\begin{document}
\title {A decomposition theorem for binary matroids with no prism minor}
\maketitle 
\begin {center}
S. R. Kingan 
\footnote{The first author is partially supported by  PSC-CUNY grant number 64181-00 42} \\     
Department of Mathematics \\
Brooklyn College, City University of New York\\
 Brooklyn, NY 11210\\
skingan@brooklyn.cuny.edu\\  
\end {center}

\begin {center}
Manoel Lemos 
\footnote{The second author is partially supported by CNPq under grant number 300242/2008-05.}\\
Departamento de Matematica \\
Universidade Federal de Pernambuco\\
Recife, Pernambuco, 50740-540, Brazil\\
manoel@dmat.ufpe.br\\  
\end {center}
\bigskip

\begin{abstract}   The prism graph is the dual of the complete graph on five vertices with an edge deleted, $K_5\backslash e$. In this paper we determine the class of binary matroids with no prism minor. The motivation for this problem is the 1963 result by Dirac where he identified the simple 3-connected graphs with no minor isomorphic to the prism graph.  We prove that besides Dirac's infinite families of  graphs and four infinite families of non-regular matroids determined by Oxley, there are only three possibilities for a matroid in this class: it is isomorphic to the dual of the generalized parallel connection of $F_7$ with itself across a triangle with an element of the triangle deleted; it's rank is bounded by 5; or it admits a non-minimal exact 3-separation induced by the 3-separation in $P_9$. Since the prism graph has rank 5, the class has to contain the binary projective geometries of rank 3 and 4, $F_7$ and $PG(3, 2)$, respectively. We show that there is just one rank 5 extremal matroid in the class. It has 17 elements and is an extension of $R_{10}$, the unique splitter for regular matroids. As a corollary, we obtain Dillon, Mayhew, and Royle's result identifying the binary internally 4-connected matroids with no prism minor [5].

\end{abstract}

\bigskip 

\section {\bf Introduction}

In a decomposition result, a more complicated matroid is broken down into simpler components. The fact that such simplifications exist is surprising and indicative of deep order in the structure of infinite classes of matroids. 
In 1980 Seymour decomposed the class of regular matroids, begining a flourishing genre of such structural results [9]. A matroid is {\it regular} if it has no minor isomorphic to the Fano matroid $F_7$ or its dual $F^*_7$. To decompose regular matroids, he developed the  Splitter Theorem, a Decomposition Theorem, and the notion of 3-sums. The Splitter Theorem describes how 3-connected matroids can be systematically built-up and the Decomposition Theorem describes the conditions under which a specific type of separation in a matroid gets carried forward to all matroids containing it. The proof of the decomposition of regular matroids consists of three main parts. The first part establishes that a 3-connected regular matroid is graphic or cographic or has a minor isomorphic to $R_{10}$ or $R_{12}$. The matroid $R_{10}$ is a splitter for regular matroids. This means no 3-connected regular matroid contains it (other than $R_{10}$ itself). So the building-up process stops at $R_{10}$. The second part establishes that $R_{12}$ has a non-minimal exact 3-separation that carries forward in all matroids containing it. The third part establishes that 3-connected regular matroids can be pieced together from graphic and co-graphic matroids using the operation of 3-sums. It is well-known that matroids that are not 3-connected can be pieced together from 3-connected matroids using the operations of 1-sum and 2-sum, so it sufficies to focus on the 3-connected members of a class.

We present the decomposition of binary matroids with no minor isomorphic to the prism graph. To decompose this class we used a strengthening of the Splitter Theorem [3] and a decomposition theorem by Mayhew, Royle, and Whittle [4]. The class of binary matroids with no prism minor is quite different from the class of regular matroids, but also similar in the sense that there are several special matroids in it and one of them has a separation that carries forward. The role of $R_{12}$ is played by the non-regular matroid $P_9$. 
The prism graph, shown in Figure 1, is the dual of the complete graph on five vertices, $K_5$ with one edge deleted. It is denoted as $(K_5\backslash e)^*$.

\begin{figure}[h]
\centering
\epsfxsize 2in \epsfbox{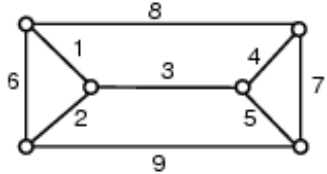}
\caption{The prism graph }
\end{figure}

\noindent A matrix representation for the prism graph is shown below.

\small
\[ 
\left[ 
\begin{array}{ccccc|cccc}
&&&& &  1&0&1&0 \\
&&I_5&&&1&0&0&1 \\
&&&&&   0&0&1&1 \\
&&&&&   0&1&1&0  \\
&&&&&   0&1&0&1
\end{array} 
\right] 
\] 
\normalsize

The origin of this excluded minor problem can be traced to 1963 when Dirac determined the extremal graphs  without two vertex disjoint cycles [1]. Excluding two vertex-disjoint cycles in a 3-connected graph is equivalent to excluding $(K_5\backslash e)^*$ as a minor.
For $r\ge 3$, let $W_r$ denote the wheel with $r$ spokes, and for $p\ge 3$,  let $K_{3,p}$ denote the complete bipartite graph with three vertices in one 
class and $p$ vertices in the other class. Let $K'_{3,p}$, $K''_{3,p}$, and $K'''_{3,p}$ denote the graphs obtained from $K_{3,p}$ by adding one, two, and three edges, respectively, joining vertices in the class containing three vertices. Dirac proved that a simple $3$-connected graph has no minor isomorphic to $(K_5\backslash e)^*$ if and only if it is isomorphic to   $W_r$ for some $r\ge 3$, $K_5$, $K_5 \backslash e$, $K_{3,p}$, $K'_{3,p} $, $K''_{3,p} $ or $K'''_{3,p}$ for some $p\ge 3$. In 1984 Robertson and Seymour published a note where they proved that a simple $3$-connected graph with at least four vertices has no minor isomorphic to $K_5\backslash e$ if and only if it is isomorphic to $(K_5 \backslash e )^*$, $K_{3,3}$,  or $W_r$ for some $r\ge 3$ [8]. In 1996 Kingan characterized the $3$-connected regular matroids with no minor isomorphic to $M^*(K_5 \backslash e)$ [2, 2.1].   

\bigskip

\noindent{\bf Theorem 1.1.} {\it M is a $3$-connected regular matroid with no minor isomorphic to $M^*(K_5 \backslash e)$ if and only if M is isomorphic to  $M(W_r)$ for some $r \ge 3$, $M(K_5)$, $M(K_5 \backslash e)$, $M^*(K_{3,3})$, $M(K_{3,p})$, $M(K'_{3,p})$, $M(K''_{3,p})$, or $M(K'''_{3,p})$ for some $p \ge 3$, or $R_{10}$.} 
\bigskip

Some matroids like $R_{10}$ play a significant role in structural results. This class contains one such significant matroid called $E_5$. It is a self-dual non-regular internally 4-connected single-element extension of $M(K_{3,3})$.  It is also  the splitter for the class  of binary matroids with no minor isomorphic to the prism graph, it dual, or the binary affine cube $AG(3,2)$ [2].  A matrix representation is shown below. 

\small
\[ 
E_5=\left[ 
\begin{array}{cccc|cccccc}
&&&&&    0&1&1&1&1 \\
&&&&&    1&0&1&1&0 \\
&&I_5&&& 1&1&0&1&1 \\
&&&&&    1&1&1&1&0 \\
&&&&&    1&1&0&0&0
\end{array} 
\right] 
\] 
\normalsize

In order to characterize the class of binary non-regular 3-connected matroids with no prism minor, we flag a particular binary non-regular 9-element rank-4 matroid known as $P_9$ and prove that besides a few exceptional matroids, all the matroids in the class have $P_9$ as a 3-decomposer. This means the matroid has a non-minimal exact $3$-separation induced by the non-minimal exact $3$-separation of $P_9$. $P_9$ is the generalized parallel connection, $P_\triangle (F_7, W_3)$, of $F_7$ and $W_3$ across a triangle, with the rim element of the triangle deleted. Note that $P_9\cong H_4$, mentioned above. A matrix representation for $P_9$ is given below. 
\small
\[ 
P_9=\left[ 
\begin{array}{cc|ccccc}
&&   0&1&1&1&1  \\
&I_4&1&0&1&1&1 \\
&&   1&1&0&1&0 \\
&&   1&1&1&1&0
\end{array} 
\right] 
\] 
\normalsize

The matroid $P_9$ appears in [5] where Oxley characterized  the $3$-connected binary non-regular matroid with no minors isomorphic to $P_9$ or $P_9^*$. Members of this class are the infinite families $Z_r$, $Z_r^*$, $Z_r\backslash b_r$ and $Z_r\backslash c_r$. The matroid $Z_r$ is a  $(2r+1)$-element rank-$r$ non-regular matroid. It can be represented by the binary matrix $[I_r| D]$ where $D$ has $r+1$ columns labeled $b_1, \dots , b_r, c_r$. The first $r$ columns in $D$ have zeros along the diagonal and ones elsewhere. The last column is all ones. Note that $Z_4\backslash c_4 \cong AG(3, 2)$ and $Z_4\backslash b_4 \cong S_8$, where $S_8$ and $AG(3, 2)$ are the two non-isomorphic single-element extensions of $F_7^*$. The next theorem appears in [6].

\bigskip
\noindent {\bf Theorem 1.2.} {\it $M$ is a $3$-connected binary non-regular matroid with no minor isomorphic to $P_9$ or $P_9^*$ if and only if $M$ is isomorphic to $F_7$, $F_7^*$, $Z_r$, $Z_r^*$, $Z_r\backslash b_r$, or  $Z_r\backslash c_r$, for some $r\ge 4$.} $\qed$
\bigskip

It is easy to show that the above infnite families do not have a prism minor. As a consequence we may conclude a binary non-regular matroid with no prism minor is either one of the infinite families mentioned in Theorems 1.2 or it has a $P_9$-minor. Like $R_{12}$, $P_9$ has a non-minimal exact 3-separation in it. However, unlike $R_{12}$, the separation in $P_9$ does not extend to all matroids containing it. Nonetheless, we are able to identify all of the exceptions.  Clearly all binary non-regular 3-connected rank 4 matroids have no  prism minor, since the prism graph has rank 5. Thus $PG(3,2)$ and all of its deletion minors have no prism minor. We prove that besides one 11-element rank-6 matroid, all the exceptions have rank at most 5.
The rank 6 exception is the dual of the generalized parallel connection of $F_7$ with itself across a triangle, with an element of the triangle deleted, denoted as $(P_\triangle (F_7, F_7)\backslash e)^*$. We are now ready to state the main results of this paper.  
\bigskip

\bigskip
\noindent{\bf Theorem 1.3.} {\it Suppose $M$ is a $3$-connected binary non-regular matroid with no $M^*(K_5 \backslash e)$-minor. Then one of the following holds:

\begin {enumerate}
\item [(i)] $M$ is isomorphic to $Z_r$, $Z_r^*$, $Z_r\backslash b_r$, or $Z_r\backslash c_r$, for some $r\ge 4$;
\item [(ii)] $P_9$ is a $3$-decomposer for $M$;
\item [(iii)] $M$ is isomorphic to $(P_\triangle (F_7, F_7)\backslash z)^*$; or
\item [(iv)] $M$ has rank at most $5$.
\end{enumerate} }
\bigskip

A detailed analysis of  rank 5 binary matroids reveals that all of them are restriction minors of one particular  17-element matroid $R_{17}$, that is an extension of $E_5$ and $R_{10}$  [2].  A matrix representation for $R_{17}$ is shown below.  

\small
\[ 
R_{17}=\left[ 
\begin{array}{c|cccccccccccc}
&    1&0&0&1&1&0&0&1&1&1&1&1 \\
&    1&1&0&0&1&1&1&0&0&1&1&1 \\
I_5& 1&1&1&0&0&0&1&1&1&0&1&1 \\
&    0&1&1&1&0&1&0&0&1&1&1&1 \\
&    0&0&1&1&1&1&1&1&0&0&1&0
\end{array}
\right] 
\] 
\normalsize

\noindent As a corollary of Theorem 1.3, we obtain the following characterization of binary matroids with no prism minor.  
\bigskip

\noindent{\bf Theorem 1.4.} {\it Suppose $M$ is a $3$-connected binary matroid with no $M^*(K_5 \backslash e)$-minor. Then  either $P_9$ is a $3$-decomposer for $M$ or $M^*$ is isomorphic to one of the following matroids:
\begin {enumerate}
\item [(i)] $M(W_r)$ for some $r\ge 3$, $M(K_{3,p})$, $M(K'_{3,p})$, $M(K''_{3,p})$ or $M(K'''_{3,p})$, for some $p\ge 3$; or
\item [(ii)] $Z_r$, $Z_r^*$, $Z_r\backslash b_r$, or $Z_r\backslash c_r$, for some $r\ge 4$;
\item [(iii)] $F_7$, $(P_\triangle (F_7, F_7)\backslash z)^*$; or
\item [(iv)] $PG(3,2)$ or $R_{17}$ or one of their $3$-connected restrictions. 
\end{enumerate} }
\bigskip

\noindent Note that, $Z_4\backslash c_4 \cong AG(3, 2)$, $Z_4\backslash b_4 \cong S_8$, $M(K_5)$, $M(K_5 \backslash e)$, $M^*(K_{3,3})$, $F_7^*$ and $P_9$ are restrictions (deletion-minors) of $PG(3, 2)$.  The matroid $(P_\triangle (F_7, F_7)\backslash z)^*$ has rank 6 and 10 elements.  Therefore, $P_\triangle (F_7, F_7)\backslash z$ is a rank 4, 10-element matroid and a restriction of $PG(3, 2)$. The matroid $R_{10}$, $P_9^*$, and $E_5$   are restrictions of $R_{17}$. So we do not have to list these matroids explicitly. 

In the next section we give the statement and proof of the decomposition theorem that forms a key component of Theorem 1.3. In Section 4 we prove Theorem 1.4. We also determine the class of binary matroids with no $M(K_5\backslash e)$-minor and the class with neither $M(K_5\backslash e)$  nor  $M^*(K_5\backslash e)$-minor. 
\bigskip

\section {\bf Proof of Theorem 1.3}
\bigskip

The matroid terminology follows Oxley [7]. We should note that the matroid corresponding to the matrix labeled A is called $M[A]$ and not just $A$. However, we refer to large numbers of matrices in this paper and with the reader's understanding treat the matrix and matroid as synonymous.

Let $M$ be a matroid and $X$ be a subset of the ground set $E$. The {\it connectivity function}  $\lambda$ is defined as $\lambda (X) = r(X) + r(E-X) - r(M)$.  Observe that $\lambda (X) = \lambda (E-X)$.  For $k\ge 1$, a partition $(A, B)$ of $E$ is called a $k$-separation if $|A|\ge k$, $|B|\ge k$, and $\lambda (A) \le k-1$.  When $\lambda (A)=k-1$, we call $(A, B)$ an {\it exact k-separation}.  When $\lambda (A)=k-1$ and $|A|=k$ or $|B|=k$ we call $(A, B)$ a {\it minimal exact k-separation}. For $n\ge 2$, we say $M$ is {\it n-connected} if $M$ has no $k$-separation for $k\le n-1$.  A matroid is {\it internally $n$-connected} if it is $n$-connected and has no non-minimal exact $n$-separations. In particular, a simple matroid is 3-connected if $\lambda (A)\ge 2$ for all partitions $(A, B)$ with $|A|\ge 3$ and $|B|\ge 3$. A 3-connected matroid is {\it internally $4$-connected}  if $\lambda (A)\ge 3$ for all partitions $(A, B)$ with $|A|\ge 4$ and $|B|\ge 4$.  For example, $E_5$ is internally 4-connected, but $P_9$ is not. In the matrix representation of $P_9$ in Section 1, it has a non-minimal exact 3-separation $(A, B)$ where $A=\{1, 2, 5, 6\}$.

Let $\mathcal M$ be a class of matroids closed under minors and isomorphisms. Let $k\ge 1$ and $N$ be a matroid belonging to $\mathcal M$ having an exact  $k$-separation $(A, B)$. Let $M\in \mathcal M$ having an $N$-minor.  We say that $N$ is a {\it $k$-decomposer} for $M$ having $(A, B)$ as an {\it inducer} provided $M$ has a $k$-separation $(X, Y)$ such that $A\subseteq X$ and $B\subseteq Y$.

It is well known that every non-regular binary matroid has a minor isomorphic to $F_7$ or $F_7^*$. Thus we may consider this our starting point for any investigation of non-regular matroids. As mentioned earlier, $AG(3, 2)$ and $S_{8}$ are the two non-isomorphic 3-connected single-element extensions of $F_7^*$. The matroid $S_8$ has two non-isomorphic 3-connected single-element extensions $P_9$ and $Z_4$. The matroid $AG(3, 2)$ has one 3-connected single-element extension, $Z_4$. As noted earlier, $P_9$ has a non-minimal exact 3-separation (and cosequently so does $P_9^*$). 

We begin by proving that $P_9$ or $P_9^*$ are 3-decomposers for a certain  class of matroids. To do so we use the following result by Mayhew, Royle, and Whittle in  [3, 2.10]. Then, we will prove the stronger statement that $P_9^*$ is not relevant and, in fact, $P_9$ is the required 3-decomposer (with one exception).  We end by showing that the rank of the exceptional matroids that do not have $P_9$ as a 3-decomposer is bounded by 5. This portion requires the Strong Splitter Theorem [3, 1.4].
\bigskip

\noindent {\bf Lemma 2.1.} {\it Suppose $\mathcal M$ is a class of matroids closed under minors and isomorphism and let $N\in \mathcal M$ be a $3$-connected matroid with $|E(N)|\ge 8$ and a non-minimal exact $3$-separation $(A, B)$ where $A$ is a $4$-element  circuit and a cocircuit. If $A$ is a circuit and cocircuit in every $3$-connected single-lement extension  and  coextension of $N$ in $\mathcal M$, then $N$ is a $3$-decomposer for every matroid in $\mathcal M$ with an $N$-minor.} $\qed$
\bigskip

The significance of the above decomposition result is that it makes it easy to determine whether or not a non-minimal exact 3-separation carries forward. Compare the criteria in this result to the original criteria in Seymour's Decomposition Theorem [8, 9.1].

\bigskip
\noindent{\bf Lemma 2.2.} {\it
Suppose $N$ is a $3$-connected proper minor of a $3$-connected matroid $M$ such that,  if  $N$ is a wheel or whirl then $M$ has no larger wheel or whirl-minor, respectively. Further, suppose $m=r(M)-r(N)$. Then there is a sequence of $3$-connected matroids $M_0,M_1,\dots,M_n$, for some integer $n\ge m$, such that
\begin{enumerate}
\item[(i)] $M_0\cong N$;
\item[(ii)] $M_n=M$;
\item[(iii)] for $k\in\{1,2,\dots,m\}$, $r(M_k)-r(M_{k-1})=1$ and $|E(M_k)-E(M_{k-1})|\le 3$; and
\item[(iv)] for $m<k\le n$, $r(M_k)=r(M)$ and $|E(M_k)-E(M_{k-1})|=1$.
\end{enumerate}
Moreover, when $|E(M_k)-E(M_{k-1})|=3$, for some $1\le k\le m$, $E(M_k)-E(M_{k-1})$ is a triad of $M_k$.}
\bigskip

The significance of the Strong Splitter Theorem is that we can obtain, up to isomorphism,  $M$ starting with $N$ and at each step doing a 3-connected single-element extension or coextension, such that at most two consecutive single-element extensions occur in the sequence (unless the rank of the matroids involved are $r$). Moreover, if two consecutive single-element extensions by elements $\{e, f\}$ are followed by a coextension by element $g$, then $\{e, f, g\}$ form a triad in the resulting matroid. This greatly reduces the computations we need to establish a bound on the rank.

\bigskip

\noindent {\bf Proof of Theorem 1.3.} Since $M$ is non-regular,  $M$ has a minor isomorphic to $F_7$ or $F_7^*$. Since $F_7$ has no binary extensions, we may assume $M$ has a minor isomorphic to $F_7^*$. The 8-element binary simple extensions of $F_7^*$ are $AG(3,2)$ and $S_8$ and the 9-element simple extensions are $P_9$ and $Z_4$. We present the proof as a series of claims.

\medskip

\noindent {\bf Claim 1.} If $M$ has no $P_9$ nor $P_9^*$-minor, then $M$ is isomorphic to $F_7$, $F_7^*$, $Z_r$,   $Z_r^*$, $Z_r\backslash c_r$, or $Z_r\backslash b_r$  for $r\ge 4$.
\bigskip

\noindent {\bf Proof.} Theorem 1.2  identifies the above families as the binary  non-regular 3-connected matroids with no $P_9$ nor $P_9^*$-minor.  To prove Theorem 1.2, Oxley proves that for $r\ge 4$, $Z_r$,   $Z_r^*$, $Z_r\backslash c_r$,  and $Z_r\backslash b_r$  have no $M(W_4)$-minor [5, Theorem 2.1]. Since $M^*(K_5\backslash e)$ and $M(K_5\backslash e)$-minor have an $M(W_4)$-minor, we may conclude that $Z_r$,   $Z_r^*$, $Z_r\backslash c_r$,  and $Z_r\backslash b_r$ have no minor isomorphic to $M^*(K_5\backslash e)$ nor $M(K_5\backslash e)$.
\medskip

Returning to the proof of the theorem, Claim 1 implies that $M$ has a minor isomorphic to $P_9$ or $P_9^*$.  Now, $P_9$ has  three simple non-isomorphic binary single-element extensions shown below. Adding column $[1 1 1 0]$ gives $D_1$,  adding columns $[1 0 0 1]$, $[0 1 0 1]$, $[0 1 1 0]$, or $[1 0 1 0]$ gives $D_2$, and adding column  $[0 0 1 1]$ gives $D_3$. This is concisely summarized in Table 1a and representative matrices for $D_1$, $D_2$, and $D_3$ are given below. Note that, Table 1a gives the extensions of $P_9$. Columns in bold are the ones used to form the matrices. The final three rank 4 matrices are $PG(3, 2)$, $PG(3, 2)\backslash e$ and $PG(3, 2)\backslash \{e, f\}$.)

\tiny
\[ 
D_1=\left[ 
\begin{array}{c|cccccc}
&   0&1&1&1&1&1 \\
I_4&1&0&1&1&1&1 \\
&   1&1&0&1&0&1 \\
&   1&1&1&1&0&0
\end{array} 
\right] 
D_2=\left[ 
\begin{array}{c|cccccc}
&   0&1&1&1&1&1 \\
I_4&1&0&1&1&1&0 \\
&   1&1&0&1&0&0 \\
&   1&1&1&1&0&1
\end{array} 
\right] 
D_3=\left[ 
\begin{array}{c|cccccc}
&   0&1&1&1&1&0 \\
I_4&1&0&1&1&1&0 \\
&   1&1&0&1&0&1 \\
&   1&1&1&1&0&1
\end{array} 
\right] 
\] 
\normalsize

\tiny
 \begin{center}
\begin{tabular}{|c|p{15em}|c|}
\hline
\bf{Matroid}& \bf{Extension Columns} & {\bf Name}   \\ \hline  
$P_9$ &  $\bf [1 1 1 0]$  & $D_1$   \\ \hline
&  $\bf [1 0 0 1]$ $[0 1 0 1]$ $[0 1 1 0]$, $[1 0 1 0]$  & $D_2$               \\ \hline
& $\bf [0 0 1 1]$    &   $D_3$          \\  \hline \hline
$D_1$ & $[0 1 0 1]$ $[0 1 1 0]$ $[1 0 0 1]$ $[1 0 1 0]$ & $X_1$ \\  \hline  
& $\bf [0 0 1 1]$ & $X_2$   \\  \hline 
$D_2$ & $\bf [1 0 1 0]$ $[1 1 1 0]$  & $X_1$   \\  \hline
& $\bf [0 0 1 1 ]$ $[0 1 0 1]$ $[0 1 1 0]$ & $X_3$   \\  \hline 
$D_3$ & $[1 1 1 0]$  & $X_2$   \\  \hline
& $ [0 1 0 1]$ $[0 1 1 0]$ $[1 0 0 1]$ $[1 0 1 0]$ & $X_3$ \\ \hline 
$X_1$ & $\bf [0 0 1 1]$ $[0 1 0 1]$ $[0 1 1 0]$  & $Y_1$   \\  \hline
& $\bf [1 1 1 0]$ & $Y_2$ \\ \hline
$X_2$ & $[0 1 0 1]$ $[0 1 1 0]$ $[1 0 0 1]$ $[1 0 1 0]$  & $Y_1$ \\ \hline 
$X_3$ & $ [0 1 0 1]$ $[0 1 1 0]$ $[1 0 1 0]$ $[1 1 1 0]$  & $Y_1$   \\  \hline 
\end{tabular}
 \end{center}
\normalsize
 \begin{center} Table 1a: Rank 4 extensions of $P_9$ \end{center} 
\medskip

\noindent $P_9$ has eight cosimple non-isomorphic single-element coextensions (see Table 1b). When coextending a rank-4 matrix the column $[0, 0, 0, 0, 1]$ is added as the fifth element and a new row is added at the bottom of the right hand side of the matrix. The coextended element is column 5. 

\tiny
\[ 
E_1=\left[ 
\begin{array}{c|ccccc}
&    0&1&1&1&1 \\
&    1&0&1&1&1 \\
I_5& 1&1&0&1&0 \\
&    1&1&1&1&0 \\
&    1&1&0&0&0
\end{array} 
\right] 
E_2=\left[ 
\begin{array}{c|ccccc}
&    0&1&1&1&1 \\
&    1&0&1&1&1 \\
I_5& 1&1&0&1&0 \\
&    1&1&1&1&0 \\
&    1&1&0&1&1
\end{array} 
\right] 
E_3=\left[ 
\begin{array}{c|ccccc}
&    0&1&1&1&1 \\
&    1&0&1&1&1 \\
I_5& 1&1&0&1&0 \\
&    1&1&1&1&0 \\
&    1&1&0&0&1
\end{array} 
\right] 
E_4=\left[ 
\begin{array}{c|ccccc}
&    0&1&1&1&1 \\
&    1&0&1&1&1 \\
I_5& 1&1&0&1&0 \\
&    1&1&1&1&0 \\
&    0&1&0&0&1
\end{array} 
\right] 
\] 
\normalsize

\tiny
\[ 
E_5=\left[ 
\begin{array}{c|ccccc}
&    0&1&1&1&1 \\
&    1&0&1&1&1 \\
I_5& 1&1&0&1&0 \\
&    1&1&1&1&0 \\
&    1&0&1&0&0
\end{array} 
\right] 
E_6=\left[ 
\begin{array}{c|ccccc}
&    0&1&1&1&1 \\
&    1&0&1&1&1 \\
I_5& 1&1&0&1&0 \\
&    1&1&1&1&0 \\
&    0&0&1&0&1
\end{array} 
\right] 
E_6^*=\left[ 
\begin{array}{c|ccccc}
&    0&1&1&1&1 \\
&    1&0&1&1&1 \\
I_5& 1&1&0&1&0 \\
&    1&1&1&1&0 \\
&    0&0&1&1&1
\end{array} 
\right] 
E_7=\left[ 
\begin{array}{c|ccccc}
&    0&1&1&1&1 \\
&    1&0&1&1&1 \\
I_5& 1&1&0&1&0 \\
&    1&1&1&1&0 \\
&    0&0&0&1&1
\end{array} 
\right] 
\] 
\normalsize

\tiny
 \begin{center}
\begin{tabular}{|p{28em}|c|}
\hline
\bf{Coextension Rows} & {\bf Name}     \\   \hline \hline
$[1 1 0 0 0]$ $[1 1 1 1 1]$    & $E_1$  \\  \hline
$[1 1 0 1 1]$ $[1 1 1 0 0]$    & $E_2$  \\  \hline
$[1 1 0 0 1]$ $[1 1 1 0 1]$    &  $E_3$  \\  \hline
$[0 1 0 0 1]$ $[0 1 0 1 0]$ $[0 1 1 0 1]$ $[0 1 1 1 0]$ $[1 0 0 0 1]$ $[1 0 0 1 0]$ $[1 0 1 0 1]$ $[1 0 1 1 0]$    &  $E_4$  \\  \hline
$[0 1 0 1 1]$ $[0 1 1 0 0]$ $[1 0 0 1 1]$ $[1 0 1 0 0]$   &  $E_5$  \\  \hline
$[0 0 1 0 1]$ $[0 0 1 1 0]$    &  $E_6$  \\  \hline
$[0 0 1 1 1]$                  &  $E_6^*$  \\  \hline
$[0 0 0 1 1]$      &  $E_7$  \\  \hline
\end{tabular}
 \end{center}
\normalsize
 \begin{center} Table 1b: Single-element coextensions of $P_9$ \end{center} 
\bigskip

\noindent {\bf Claim 2.}  If $M$ has a $P_9$-minor, but no $D_2$, $D_2^*$, $E_4$, or $E_5$-minor, then $P_9$ or $P_9^*$ is a $3$-decomposer for $M$.  
\bigskip

\noindent {\bf Proof.} As mentioned earlier,  $P_9$ has  a non-minimal exact $3$-separation $(A, B)$ where $A=\{1, 2, 5, 6\}$ is both a circuit and a cocircuit.   It is easy to check that the set $A=\{1, 2, 5, 6\}$ is both a circuit and a cocircuit in $D_1$ and $D_3$, whereas $D_2$ is internally 4-connected. 
Next, the set $A=\{1, 2, 5, 6\}$ corresponds to $A'=\{1, 2, 6, 7\}$ in the coextension since the fifth column is the coextended element. It can be checked that $\{1, 2, 6, 7\}$ is both a circuit and a cocircuit in $E_1$, $E_2$, $E_3$, $E_6$, $E_6^*$, and $E_7$.  Since $E_4$ and $E_5$ are self-dual, the claim follows from Lemma 2.1.   
\bigskip

\noindent {\bf Claim 3.} If $M$ has  no $D_2$, $D_2^*$, $E_4$ or $E_5$-minor,  then one of the following hold:
\begin {enumerate}
\item [(i)] $M$ is isomorphic to $F_7$, $F_7^*$,  $Z_r$, $Z_r^*$, $Z_r\backslash b_r$, or $Z_r\backslash c_r$, for some $r\ge 4$; or
\item [(ii)]$P_9$ or $P_9^*$ is a $3$-decomposer for $M$ 
\end {enumerate}
   
\bigskip

\noindent {\bf Proof.}  If $M$ has no $P_9$ or $P_9^*$-minor, then (i) follows from Claim 1. So, we may assume that $M$ has a $P_9$ or $P_9^*$ as a minor. If $M$ has a $P_9$-minor, then $P_9$ is a 3-decomposer for $M$. Otherwise, $M$ has a $P_9^*$-minor and by duality $P_9^*$ is a 3-decomposer for $M$.
\bigskip

\noindent {\bf Claim 4.} If  $M$ has an $E_5$-minor, then either $M\cong E_5$ or $M$ has a $D_2$ or $D_2^*$-minor. 
\bigskip

\noindent {\bf Proof.} The matroid $E_5$ is self-dual and has seven non-isomorphic binary 3-connected single-element extensions, shown in Table 2a. All of them have a minor isomorphic to $D_2$ and the claim follows.

\tiny
 \begin{center}
\begin{tabular}{|p{15em}|c|c|c|}
\hline
 \bf{Extension Columns} & {\bf Name} & {\bf Contraction-minor}  & {\bf Deletion-minor}\\      \hline \hline
  $[0 0 1 0 1]$ $[0 0 1 1 0]$ $[0 1 0 1 1]$ $[0 1 1 0 0]$  & $A$              & $D_2$ &  $E_5$ $E_6^*$, $E_7$, $K'_{3,3}$     \\  \hline

 $[1 0 0 1 1]$ & $B$       &    $D_2$ & $E_5$, $K'_{3,3}$, $R_{10}$ \\  \hline

 $[1 1 0 0 1]$ $[1 1 1 0 1]$   &   $C$            &  $D_2$ & $E_5$, $E_6^*$, $E_7$, $E_3$    \\  \hline

 $[0 0 0 1 1 ]$ $[0 0 1 1 1]$ $[0 1 0 0 1]$ $[0 1 1 0 1]$     &   & $D_2$ & $E_4$    \\  \hline

 $0 1 0 1 0]$ $[0 1 1 1 0]$  &    & $D_2$ & $E_4$    \\  \hline

 $[1 0 0 0 1]$ $[1 0 0 1 0]$  $[1 1 0 1 1]$ $[1 1 1 0 0]$  &                &  $D_2$ & $E_4$    \\  \hline

 $[1 0 1 0 1]$ $[1 0 1 1 0]$ $[1 1 0 0 0]$ $[1 1 1 1 1]$                &            &  $D_2$ & $E_4$    \\  \hline
\end{tabular}
 \end{center}
\normalsize
 \begin{center} Table 2a: Single-element extensions of $E_5$ \end{center} 
\bigskip

Returning to the proof of the theorem, we will determine which of the extensions and coextensions of $P_9$ have minors isomorphic to $M^*(K_5\backslash e)$ or $M(K_5\backslash e)$. A matrix representation for the graph $K_5\backslash e$ is given below.

\tiny
\[ 
K_5\backslash e=\left[ 
\begin{array}{c|ccccc}
&   1&0&0&1&1 \\
I_4&1&1&0&0&0 \\
&   0&1&1&0&1 \\
&   0&0&1&1&0
\end{array} 
\right] 
\]
\normalsize
\noindent It has three binary non-isomorphic single-element extensions, $K_5$, $D_2$, and $D_3$ (see Table 3a).  

\tiny
\[ 
K_5=\left[ 
\begin{array}{c|cccccc}
&   1&0&0&1&1&0 \\
I_4&1&1&0&0&0&1 \\
&   0&1&1&0&1&0 \\
&   0&0&1&1&0&1
\end{array} 
\right] 
D_2=\left[ 
\begin{array}{c|cccccc}
&   1&0&0&1&1&0 \\
I_4&1&1&0&0&0&1 \\
&   0&1&1&0&1&1 \\
&   0&0&1&1&0&1
\end{array} 
\right] 
D_3=\left[ 
\begin{array}{c|cccccc}
&   1&0&0&1&1&1 \\
I_4&1&1&0&0&0&0 \\
&   0&1&1&0&1&1 \\
&   0&0&1&1&0&1
\end{array} 
\right] 
\] 
\normalsize

\tiny
 \begin{center}
\begin{tabular}{|p{15em}|c|}
\hline
\bf{Extension Columns} & {\bf Name}  \\      \hline \hline
$[0 1 0 1]$   & $K_5$   \\  \hline
$[0 1 1 1]$ $[1 1 0 1]$ $[1 1 1 1]$  & $D_2$               \\  \hline
$[1 0 1 1]$ $[1 1 1 0]$    &   $D_3$          \\  \hline
\end{tabular}
 \end{center}
\normalsize
 \begin{center} Table 3a: Single-element extensions of $M(K_5\backslash e)$ \end{center} 
\bigskip

\noindent The isomorphisms from these representations of $D_2$ and $D_3$ to the previous representations are, respectively, $$\{1, 2, 3, 4, 5, 6, 7, 8, 9, 10\} \rightarrow \{3, 4, 5, 8, 2, 10, 9, 6, 7, 1\}$$ and $$\{1, 2, 3, 4, 5, 6, 7, 8, 9, 10\} \rightarrow \{8, 7, 9, 1, 3, 4, 2, 5, 10, 6 \}.$$  

\noindent Using the matrix representation of $M^*(K_5\backslash e)$ in the previous section we see that up to isomorphism the 3-connected binary single-element extensions are as follows:

\tiny
\[ 
G=\left[ 
\begin{array}{c|ccccc}
&    1&0&1&0&0 \\
&    1&0&0&1&0 \\
I_5& 0&0&1&1&1 \\
&    0&1&1&0&0 \\
&    0&1&0&1&1
\end{array} 
\right] 
(K_{3,3}')^*=\left[ 
\begin{array}{c|ccccc}
&    1&0&1&0&0 \\
&    1&0&0&1&1 \\
I_5& 0&0&1&1&0 \\
&    0&1&1&0&0 \\
&    0&1&0&1&1
\end{array} 
\right] 
\] 
\normalsize

\tiny
\[ 
E_4=\left[ 
\begin{array}{c|ccccc}
&    1&0&1&0&0 \\
&    1&0&0&1&1 \\
I_5& 0&0&1&1&0 \\
&    0&1&1&0&1 \\
&    0&1&0&1&0
\end{array} 
\right] 
E_6=\left[ 
\begin{array}{c|ccccc}
&    1&0&1&0&0 \\
&    1&0&0&1&0 \\
I_5& 0&0&1&1&1 \\
&    0&1&1&0&1 \\
&    0&1&0&1&1
\end{array} 
\right] 
E_7^*=\left[ 
\begin{array}{c|ccccc}
&    1&0&1&0&1 \\
&    1&0&0&1&1 \\
I_5& 0&0&1&1&1 \\
&    0&1&1&0&1 \\
&    0&1&0&1&1
\end{array} 
\right] 
\] 
\normalsize

\tiny
 \begin{center}
\begin{tabular}{|p{20em}|c|}
\hline
\bf{Extension Columns} & {\bf Name}  \\      \hline \hline
$[0 0 1 0 1]$ $[0 0 1 1 0]$ $[0 1 1 0 0]$ $[0 1 1 1 0]$ $[1 0 1 0 0]$ $[1 0 1 0 1]$   & $(K_5\backslash e)^*+edge$   \\  \hline
$[0 1 0 0 1]$ $[1 0 0 1 0]$ $[1 1 0 1 1]$    &   $(K'_{3,3})^*$          \\  \hline
$[0 1 0 1 0]$ $[0 1 0 1 1]$ $[1 0 0 0 1]$ $[1 0 0 1 1]$ $[1 1 0 0 1]$ $[1 1 0 1 0]$ & $E_4$               \\  \hline
$[0 0 1 1 1]$ $[0 1 1 1 1]$ $[1 0 1 1 1]$ $[1 1 1 0 0]$ $[1 1 1 0 1]$ $[1 1 1 1 0]$  & $E_6$               \\  \hline
$[1 1 1 1 1]$  & $E_7^*$              \\  \hline
\end{tabular}
 \end{center}
\normalsize
 \begin{center} Table 3b: Single-element extensions of $M(K_5\backslash e)$ \end{center} 
\bigskip

\noindent The isomorphisms from these representations of $E_4$, $E_6$ and $E_7^*$ to the representations of $E_4$, $E_6$ and dual of $E_7$  are, respectively, $$\{1, 2, 3, 4, 5, 6, 7, 8, 9, 10\} \rightarrow \{3, 9, 2, 6, 7, 8, 10, 4, 5, 1\}$$  $$\{1, 2, 3, 4, 5, 6, 7, 8, 9, 10\} \rightarrow \{4, 10, 1, 7, 2, 8, 9, 3, 5, 6\}$$ and $$\{1, 2, 3, 4, 5, 6, 7, 8, 9, 10\} \rightarrow \{3, 4, 1, 9, 10, 8, 5, 2, 6, 7 \}  $$  

\noindent Thus we conclude that $M^*(K_5\backslash e)$ has five non-isomorphic binary $3$-connected single element extensions, the graph $G$ obtained by adding an edge to $(K_5\backslash e)^*$, the cograph $M^*(K_{3, 3})$, and three binary non-regular matroids $E_4$, $E_6$, and $E_7^*$.   

Since $D_2^*$, $E_4$, and $E_7^*$ have an $M^*(K_5\backslash e)$-minor, it follows from Claim 3 that, if $M$ has a $P_9$ or $P_9^*$-minor and rank at least 5, then either they are 3-decomposers for $M$ or $M$ has an $E_5$ or $D_2$-minor.

Now, Table 2 implies that all the extensions of $E_5$ have a $D_2$-minor, which in turn has a $M(K_5\backslash e)$-minor.  Further, all single-element extensions, except $A$, $B$, and $C$, have an $E_4$-minor, which has an $M^*(K_5\backslash e)$-minor. 
Moreover, since the extensions of $E_5$ have a $D_2$-minor and $E_5$ is self-dual, all the coextensions have a  $D_2^*$-minor (which has an $M^*(K_5\backslash e)$-minor). 
\bigskip

Matrix representations for $A$, $B$, and $C$ are given below. In Claim 5 we show that the coextensions of $A$, $B$, and $C$ have an  $M^*(K_5\backslash e)$-minor. 

\tiny
\[ 
A=\left[ 
\begin{array}{c|cccccc}
&    0&1&1&1&1&0 \\
&    1&0&1&1&0&0 \\
I_5& 1&1&0&1&1&1 \\
&    1&1&1&1&0&0 \\
&    1&1&0&0&0&1
\end{array} 
\right] 
B=\left[ 
\begin{array}{c|ccccccc}
&    0&1&1&1&1&1 \\
&    1&0&1&1&0&0 \\
I_5& 1&1&0&1&1&0 \\
&    1&1&1&1&0&1 \\
&    1&1&0&0&0&1
\end{array} 
\right] 
C=\left[ 
\begin{array}{c|ccccccc}
&    0&1&1&1&1&1 \\
&    1&0&1&1&0&1 \\
I_5& 1&1&0&1&1&0 \\
&    1&1&1&1&0&0 \\
&    1&1&0&0&0&1
\end{array} 
\right] 
\] 
\normalsize

\bigskip

\noindent {\bf Claim 5.} If $M$ is a coextension of $A$, $B$, or $C$, then $M$ has an $M^*(K_5\backslash e)$-minor.
\bigskip

\noindent {\bf Proof.}  Since $E_5$ is self-dual and every extension has a $D_2$-minor, it follows that every coextension has a $D_2^*$-minor, and consequently an $M^*(K_5\backslash e)$-minor. Suppose $M$ is a coextension of $A$, $B$, $C$. Then a partial matrix representation for $M$ is shown in Figure 2. 

\begin{figure}[h]
\centering
\epsfxsize 2in \epsfbox{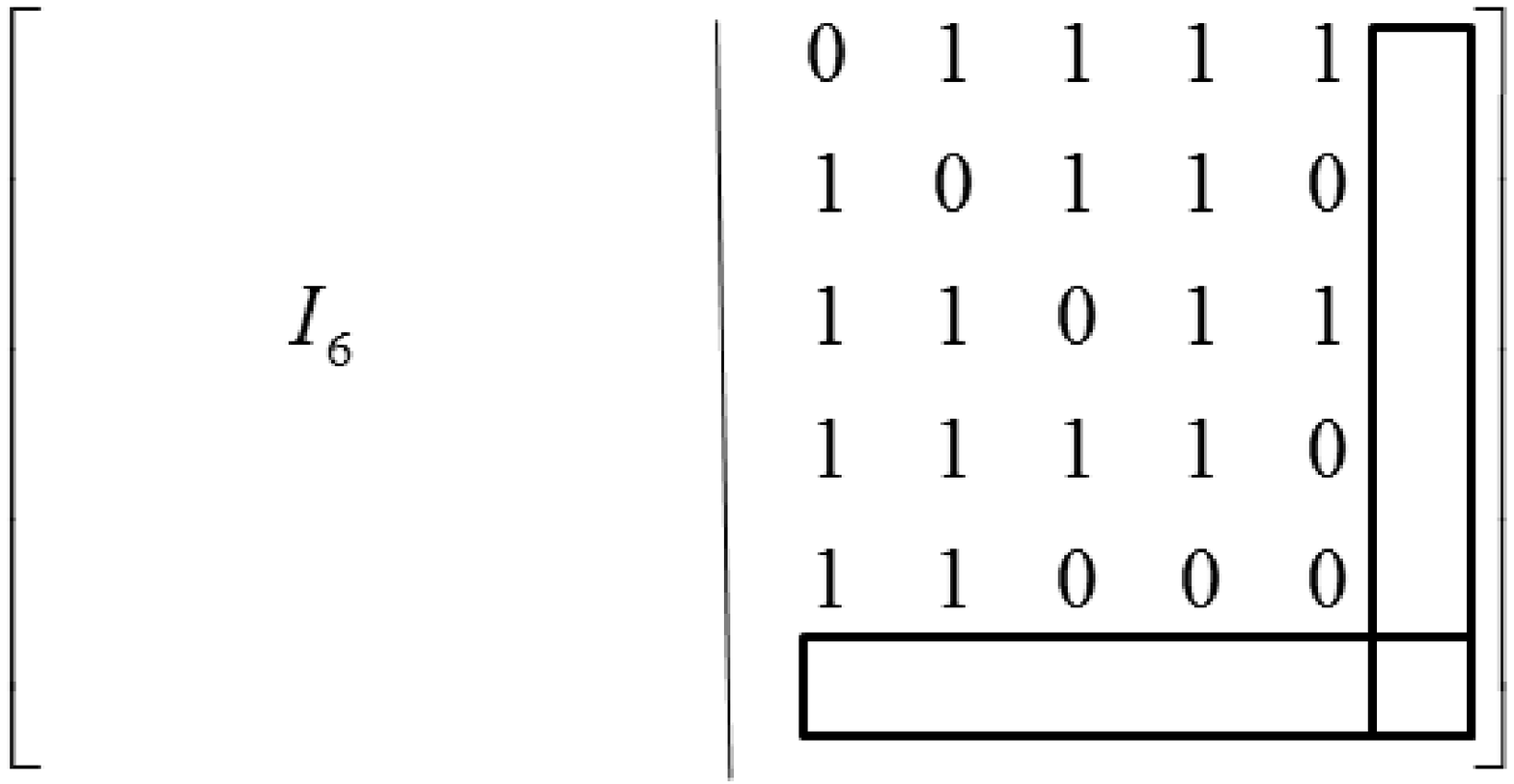}
\caption{Structure of a coextension of $A$, $B$, $C$ }
\end{figure}

There are three types of rows that may be inserted into the last row on the right-hand side of the matrix in Figure 2.

\begin {enumerate}
\item[(i)] rows that can be added to $E_5$ to obtain a coextension with no $M^*(K_5\backslash e)$-minor with a 0 or 1 as the last entry; 
\item[(ii)] the identity rows with a 1 in the last position; 
 \item[(iii)] and the rows ``in-series" to the right-hand side of matrices $A$, $B$, $C$ with the last entry reversed. 
\end{enumerate}

\noindent There are no Type I rows.
Type II rows are $[1 0 0 0 0 1]$, $[0 1 0 0 0 1]$, $[0 0 1 0 0 1]$, $[0 0 0 1 0 1]$, and $[0 0 0 0 1 1]$. 
Type III rows are specific to the matrices $A$, $B$, $C$. For matrix $A$ they are $[0 1 1 1 1 1]$, $[1 0 1 1 0 1]$, $[1 1 0 1 1 0]$, $[1 1 1 1 0 1]$, $[1 1 0 0 0 0]$. For matrix $B$ they are $[0 1 1 1 1 0]$, $[1 0 1 1 0 1]$, $[1 1 0 1 1 1]$, $[1 1 1 1 0 0]$, and $[1 1 0 0 0 0]$. For $C$ they are $[0 1 1 1 1 0]$, $[1 0 1 1 0 0]$, $[1 1 0 1 1 1]$, $[1 1 1 1 0 1]$, and $1 1 0 0 0 0]$.
Thus we see that only ten rows may be added. Table A1 in the Appendix shows that most of these rows  result in matroids that are isomorphic to matroids with an $M^*(K_5\backslash e)$-minor. Only two coextensions must be specifically checked for an $M^*(K_5\backslash e)$-minor: $(C, coextn 9)$ and $(C, coextn 10)$. 
Observe that, 
$(C, coextn 9)/12\backslash 1 \cong E_4$, and
$(C, coextn 10)/12\backslash 10 \cong E_4$. 
Since $E_4$ has an $M^*(K_5\backslash e)$-minor, we may conclude these matroids have it too.
\bigskip

\noindent {\bf Claim 6:} {\it If $M$ has a $P_9^*$-minor, but no $D_2$, $D_2^*$ $E_4$ or $E_5$-minor, then either  $P_9$ is a $3$-decomposer for $M$ or $M\cong D_1^*$.}

\bigskip

\noindent {\bf Proof.} Suppose $M$ is an extension of $P_9^*$. The extensions of $P_9^*$ are the duals of the extensions of $P_9$. Thus, from Table 1b they are $E_1$, $E_2$, $E_3$, $E_4$, $E_5$, $E_6$, $E_6^*$ and $E_7^*$. All of these matroids except $E_7^*$ have a $P_9$-minor since $E_1$, $E_2$, $E_3$,  $E_4$, and $E_5$ are self-dual and $E_6$ and $E_6^*$ are both coextensions of $P_9$. Since $E_7^*$ has an $M^*(K_5\backslash e)$-minor, we may conclude that $M$ cannot have an $E_7^*$-minor. So $P_9$ is a decomposer for $M$.     

Suppose $M$ is a coextension of $P_9^*$. Then since $D_1$, $D_2$, $D_3$ are extensions of $P_9$, $D_1^*$, $D_2^*$, and $D_3^*$ are coextensions of $P_9^*$. Of these, $M$ cannot be $D_2^*$ by hypothesis and $D_3^*$ since it has a $M(K_5\backslash e)$-minor. Therefore, we may suppose $M$ has a minor isomorphic to $D_1^*$. 

If $M\cong D_1^*$, then we are done. From Table 1a we see that $D_1$ is formed by adding just one column to $P_9$ ($[1 1 1 0]$), so any extension of $D_1$ will have a $D_2$ or $D_3$-minor. Thus any coextension of $D_1^*$ will have a $D_2^*$ or $D_3^*$-minor, which have an $M^*(K_5\backslash e)$-minor. 

The extensions of $D_1^*$ are the duals of the coextensions of $D_1$. Observe from Table 4 that all except the second coextension have a $P_9^*$-minor. Since the second coextension has an $E_7$-minor, its dual has an $E_7^*$-minor. Thus we may conclude that $M\cong D_1^*$. Note that  $D_1^* \cong (P_\triangle (F_7, F_7)\backslash z)$.
\bigskip

\tiny
 \begin{center}
\begin{tabular}{|c|p{25em}|c|c|}
\hline
Matroid &\bf{Coextension Rows} & {\bf Name} & Relevant minors  \\      \hline \hline

$D_1$ & $[0 0 0 0 1 1]$ $[0 0 0 1 0 1]$ $[0 0 1 0 1 0]$ $[0 0 1 1 0 0]$ $[0 1 0 0 1 0]$ $[0 1 0 1 0 0]$ $[0 1 1 0 1 1]$ $[0 1 1 1 0 1]$ $[1 0 0 0 1 0]$ $[1 0 0 1 0 0]$ $[1 0 1 0 1 1]$ $[1 0 1 1 0 1]$ $[1 1 0 0 0 1]$ $[1 1 0 1 1 1]$ $[1 1 1 0 0 0]$ $[1 1 1 1 1 0]$   &                & $E_1$, $E_2$, $E_3$, $E_4$, $E_6^*$     \\  \hline

& $[0 0 0 1 1 0]$     &                &   $E_7$   \\  \hline

& $[0 0 0 1 1 1]$ $[0 0 1 1 1 0]$ $[0 1 0 1 1 0]$ $[0 1 1 0 0 1]$ $[1 0 0 1 1 0]$ $[1 0 1 0 0 1]$ $[1 1 0 0 1 1]$ $[1 1 1 0 1 0]$   &     & $E_3$, $E_5$, $E_6^*$, $E_7$  \\  \hline

& $[0 0 1 0 0 1]$ $[0 0 1 1 1 1]$  &             &  $E_2$, $E_6^*$    \\  \hline

& $[0 0 1 0 1 1]$ $[0 0 1 1 0 1]$  &             &  $E_2$, $E_3$ $E_5$    \\  \hline

& $[0 1 0 0 0 1]$ $[0 1 0 0 1 1]$  $[0 1 0 1 0 1]$ $[0 1 0 1 1 1]$ $[0 1 1 0 0 0]$ $[0 1 1 0 1 0]$ $[0 1 1 1 0 0]$ $[0 1 1 1 1 0]$ $[1 0 0 0  1]$ $[1 0 0 0 1 1]$ $[1 0 0 1 0 1]$ $[1 0 0 1 1 1]$ $[1 0 1 0 0 0]$ $[1 0 1 0 1 0]$ $[1 0 1 1 0 0 ]$ $[1 0 1 1 1 0]$ &   & $E_2$, $E_6^*$   \\  \hline

& $[1 1 0 0 0 0]$ $[1 1 0 1 0 0]$  $[1 1 1 1 0 1]$ $[1 1 1 1 1 1]$   &   &      $E_1$, $E_2$ \\  \hline

& $[1 1 0 0 1 0]$ $[1 1 0 1 1 0]$  $[1 1 1 0 0 1]$ $[1 1 1 0 1 1]$    &   &  $E_2$, $E_3$  \\  \hline \hline

$D_2$ & $[0 0 0 0 1 1]$ $[0 0 0 1 0 1]$ $[0 0 0 1 1 0]$ $[0 0 1 1 1 1]$  $[1 0 0 1 1 1]$ $[1 0 1 0 0 0]$   & $A$              & $E_5$ $E_6^*$ $E_7$ $K'_{3,3}$   \\  \hline

& $[0 1 1 0 0 1]$     & $B$              &  $E_5$, $K'_{3, 3}$,$R_{10}$   \\  \hline

& $[0 1 0 1 1 1]$ $[1 1 0 0 1 1]$ $[1 1 1 0 1 0]$    &   $C$            & $E_3$ $E_5$ $E_6^*$, $E_7$    \\  \hline

& $[0 0 0 1 1 1]$  &   $Z$            & $E_7$, $R_{10}$    \\  \hline

& $[0 0 1 0 0 1]$ $[1 0 0 1 0 0]$  $[1 0 1 1 0 1]$  &   & $E_4$    \\  \hline

& $[0 0 1 0 1 0]$ $[0 0 1 1 0 0]$  $[1 0 0 0 0 1]$ $[1 0 0 0 1 0]$ $[1 0 1 0 1 1]$ $[1 0 1 1 1 0]$ &   & $E_4$    \\  \hline

& $[0 0 1 0 1 1]$ $[0 0 1 1 0 1]$  $[1 0 0 1 0 1]$ $[1 0 0 1 1 0]$ $[1 0 1  0 0 1]$ $[1 0 1 1 0 0]$  &   & $E_4$    \\  \hline

& $[0 0 1 1 1 0]$ $[1 0 0 0 1 1]$  $[1 0 1 0 1 0]$  &   & $E_4$    \\  \hline

& $[0 1 0 0 0 1]$ $[0 1 1 0 0 0]$  $[0 1 1 0 1 1]$ $[0 1 1 1 0 1]$ $[1 1 0 1 1 0]$ $[1 1 1 0 0 1]$ &   & $E_4$    \\  \hline

& $[0 1 0 0 1 0]$ $[0 1 0 1 0 0]$  $[1 1 0 0 0 0]$ $[1 1 0 1 0 1]$ $[1 1 1 1 0 0]$ $[1 1 1 1 1 1]$ &   & $E_4$    \\  \hline

& $[0 1 0 0 1 1]$ $[0 1 0 1 0 1]$  $[1 1 0 0 1 0]$ $[1 1 0 1 1 1]$ $[1 1 1 0 0 0]$ $[1 1 1 0 1 1]$ &   & $E_4$    \\  \hline

& $[0 1 0 1 1 0]$ $[0 1 1 0 1 0]$  $[0 1 1 1 0 0]$ $[0 1 1 1 1 1]$ $[1 1 0 0 0 1]$ $[1 1 1 1 1 0]$ &   & $E_4$    \\  \hline

\end{tabular}
 \end{center}
\normalsize
 \begin{center} Table 4: Single-element coextensions of $D_1$ and $D_2$ \end{center} 
\bigskip

\noindent {\bf Claim 7.} If $M$ is a coextension of $D_2$, then either $M$ has an $M^*(K_5\backslash e)$-minor or $M$ is isomorphic to $A$, $B$, $C$, or $Z$.
\bigskip

\noindent {\bf Proof.} Table 4 verifies that all the single-element coextensions of $D_2$ except for $A$, $B$, $C$, and $Z$ have an $E_4$-minor. Further, observe that the choice of rows for $Z$ is just one, so its coextensions will have a minor isomorphic to one of the other matroids. Claim 5 already established that coextensions of $A$, $B$, and $C$ have an $M^*(K_5\backslash e)$-minor. So $Z$ does not give rise to new coextensions. 
\bigskip

Claims 5, 6, and 7 and the fact that $E_7^*$ has an $M^*(K_5\backslash e)$-minor imply that there are only three possibilities for $M$: $P_9$ is a 3-decomposer for $M$; $M\cong D_1^*$; or $M$ has a minor isomorphic to $E_5$ or $D_2$.  

Returning to the proof of the theorem,  we must show that if $M$ has an $E_5$ or $D_2$-minor and no $M^*(K_5\backslash e)$-minor, then the rank of $M$ is bounded above by 5.
To do this, let us begin by computing the single-element extensions of $A$, $B$, $C$ and $Z$ with no $M^*(K_5\backslash e)$-minor. From Table 2, we may conclude that the only columns that can be added to $E_5$ to obtain a matroid with no $M^*(K_5\backslash e)$-minor are  $[0 0 1 0 1]$, $[0 0 1 1 0]$, $[0 1 0 1 1]$, $[0 1 1 0 0]$ $[1 0 0 1 1]$, $[1 1 0 0 1]$, $[1 1 1 0 1]$. Adding these columns gives us four non-isomorphic single-element extensions of $A$, $B$, and $C$. They are $D$, $E$, $F$, and $G$ shown below (all the extensions of $E_5$ are shown in Table 5). 

\tiny
\[ 
D=\left[ 
\begin{array}{c|ccccccc}
&    0&1&1&1&1&0&0 \\
&    1&0&1&1&0&0&0 \\
I_5& 1&1&0&1&1&1&1 \\
&    1&1&1&1&0&0&1 \\
&    1&1&0&0&0&1&0
\end{array} 
\right] 
E=\left[ 
\begin{array}{c|ccccccc}
&    0&1&1&1&1&0&0 \\
&    1&0&1&1&0&0&1 \\
I_5& 1&1&0&1&1&1&0 \\
&    1&1&1&1&0&0&1 \\
&    1&1&0&0&0&1&1
\end{array} 
\right] 
\] 
\normalsize

\tiny
\[ 
F=\left[ 
\begin{array}{c|ccccccc}
&    0&1&1&1&1&0&1 \\
&    1&0&1&1&0&0&1 \\
I_5& 1&1&0&1&1&1&0 \\
&    1&1&1&1&0&0&0 \\
&    1&1&0&0&0&1&1
\end{array} 
\right] 
G=\left[ 
\begin{array}{c|ccccccc}
&    0&1&1&1&1&0&1 \\
&    1&0&1&1&0&0&1 \\
I_5& 1&1&0&1&1&1&1 \\
&    1&1&1&1&0&0&0 \\
&    1&1&0&0&0&1&1
\end{array} 
\right] 
\] 
\normalsize

Suppose $M$ is a coextension of $D$, $E$, $F$, or $G$. Then the structure of $M$ is shown in Figure 3. 

\begin{figure}[h]
\centering
\epsfxsize 2in \epsfbox{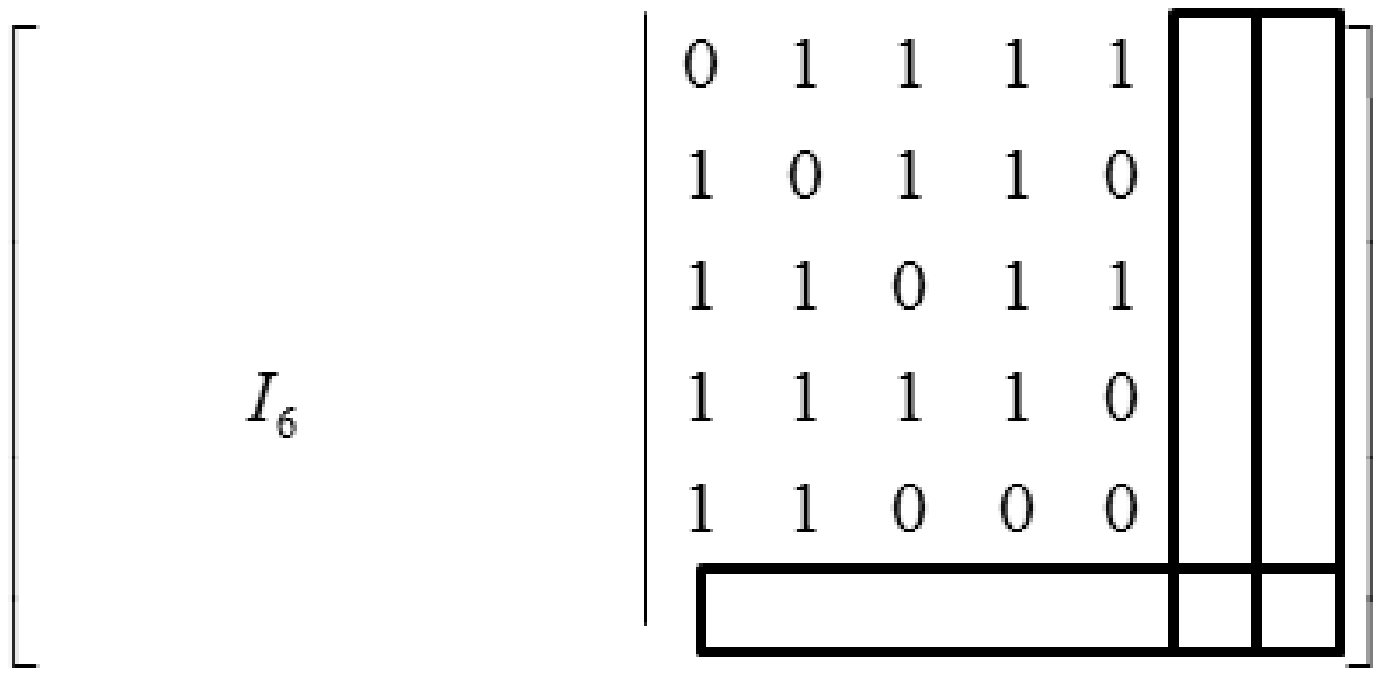}
\caption{Structure of a coextension of $D$, $E$, $F$, $G$ }
\end{figure}

\noindent Recall from the proof of Claim 5, that there are no Type I rows to add. Adding a Type II or III row (with the exception of $[0 0 0 0 0 1 1]$) causes $M\backslash 12$ to be 3-connected (and there are no such matroids).  So the  only coextension we must check is the one formed by adding row $[0 0 0 0 0 1 1]$. That is the coextensions in which $\{6, 11, 12\}$ is a triad. Let $D'$, $E'$, $F'$, and $G'$ be the coextension of $D$, $E$, $F$, and $G$, respectively, obtained by coextending by row $[0 0 0 0 0 1 1]$. Then in each case we can find an $E_4$ minor. In particular, 
$D'/1\backslash \{3, 11\} \cong E_4$,  
$E'/1\backslash \{7, 11\} \cong E_4$, 
$F'/1\backslash \{7, 11\} \cong E_4$, and
$G'/1\backslash \{7, 11\} \cong E_4$.  
Finally, observe that if $M$ is an extension of $E_5$ of size $k\ge 13$, then for some $e\in \{11, \dots , k\}$, $M\backslash e$ is 3-connected.

Next, suppose $M$ is a coextension of $Z$. Observe from Table 5 that $Z$ is an extension of $E_7$ and $R_{10}$. In the representation of $R_{17}$ given in the introduction, $R_{10}$ is isomorphic to the first ten columns and $Z$ is isomorphic to the first eleven columns. Let us take that as a representation of $Z$.

\small
\[ 
Z=\left[ 
\begin{array}{c|cccccc}
&    1&0&0&1&1&0 \\
&    1&1&0&0&1&1\\
I_5& 1&1&1&0&0&0 \\
&    0&1&1&1&0&1 \\
&    0&0&1&1&1&1
\end{array} 
\right] 
\] 
\normalsize

\noindent $R_{10}$ has two non-isomorphic binary 3-connected single-element extensions, $Z$ and $B$. $Z$ is obtained by adding any one of the columns $[0 1 0 1 1]$, $[0 1 1 0 1]$ $[1 0 1 0 1]$ $[1 0 1 1 0]$ $[1 1 0 1 0]$ or $[1 1 1 1 1]$ and $A$ is obtained by adding any one of the remaining columns (see [8] for details). Observe from the representation of $R_{17}$ in the Introduction that adding all of the above six columns to $Z$ gives us $R_{17}\backslash e$. Adding one additional column gives us $R_{17}$.

Let $M$ be a coextension of $Z$. Observe that $R_{10}$ has two single-element coextensions, $A^*$ and $Z^*$ both of which have a $D_2^*$-minor (and consequently an $M^*(K_5\backslash e)$-minor. Thus every coextension of $R_{10}$ has an $M^*(K_5\backslash e)$-minor). As  before, there are three types of rows that may be added to $M$. There are no Type I rows. Type II rows are $[1 0 0 0 0 1]$, $[0 1 0 0 0 1]$, $[0 0 1 0 0 1]$, $[0 0 0 1 0 1]$ and $[0 0 0 0 1 1]$and  Type III rows are $[1 0 0 1 1 1]$, $[1 1 0 0 1 0]$, $[1 1 1 0 0 1]$, $[0 1 1 1 0 0]$ and $[0 0 1 1 1 0]$. Observe that, adding any of the above ten rows to $Z$ gives an isomorphic matroid (Appendix Table A1). Without loss of generality let $M$ be obtained from $Z$ by adding row $[0 0 0 0 1 1]$. Then, $M/1 \backslash 7 \cong E_4$. Therefore, every coextension of $Z$ has an $M^*(K_5\backslash e)$-minor. 

It is easy to check that $Z$ has three non-isomorphic single-element extensions, {\it namely}, $D$ and $F$ mentioned above, and $Y$ shown below (Table 5).
\small
\[
Y=\left[ 
\begin{array}{c|ccccccc}
&    1&0&0&1&1&0&0 \\
&    1&1&0&0&1&1&1\\
I_5& 1&1&1&0&0&0&1 \\
&    0&1&1&1&0&1&0 \\
&    0&0&1&1&1&1&1
\end{array} 
\right] 
\] 
\normalsize
\noindent Let $Y'$ be coextension of $Y$ formed by adding row $[0 0 0 0 0 1 1]$. Then, $Y'/1\backslash \{2, 7\} \cong E_4$
Thus, we may conclude that if $M$ has an $R_{10}$-minor and no minor isomorphic to $M^*(K_5\backslash e)$-minor, then $M$ has rank at most 5.

Next, suppose $M$ has a $D_2$-minor and no minor isomorphic to $M^*(K_5\backslash e)$-minor. Then, by Claim 7  the only coextensions of $D_2$ with no $M^*(K_5\backslash e)$-minor are $A$, $B$, $C$, and $Z$. From  Table 1a we see that $D_2$ has two single-element extensions $X_1$ and $X_3$ shown below:

 \small
\[ 
X_1=\left[ 
\begin{array}{c|ccccccc}
&   0&1&1&1&1&1&1\\
I_4&1&0&1&1&1&0&0 \\
&   1&1&0&1&0&0&1 \\
&   1&1&1&1&0&1&0
\end{array} 
\right] 
X_3=\left[ 
\begin{array}{c|ccccccc}
&   0&1&1&1&1&1&0\\
I_4&1&0&1&1&1&0&0 \\
&   1&1&0&1&0&0&1 \\
&   1&1&1&1&0&1&1
\end{array} 
\right] 
\] 
\normalsize

There are three types of rows that may be added to $X_1$ and $X_3$. 

\begin {enumerate}
\item[(i)] the  rows that can be added to $D_2$ to obtain a coextension with no $M^*(K_5\backslash e)$-minor with a 0 or 1 in the last entry.  (These are the rows corresponding to $A$, $B$, $C$, $Z$ in Table 4.)
\item[(ii)] the identity rows with a 1 in the last position; 
 \item[(iii)] and the rows ``in-series" to the right-hand side of matrices $X_1$ and $X_3$ with the last entry reversed. 
\end{enumerate}

\noindent If $M$ is the coextension obtained by adding the first type of row, then $M\backslash 12$ is isomorphic to $A$, $B$, $C$, or $Z$.  Thus $M$ is either $D$, $E$, $F$, $G$ or $Y$. 
Type II rows are $[1 0 0 0 0 0 1]$, $[0 1 0 0 0 0 1]$, $[0 0 1 0 0 0 1]$, $[0 0 0 1 0 0 1]$, $[0 0 0 0 1 0 1]$, $[0 0 0 0 0  1 1]$. 
Type III rows for $X_1$ are $[0 1 1 1 1 1]$, $[1 0 1 1 0 1]$, $[1 1 0 1 1 0]$, $[1 1 1 1 0 1]$, $[1 1 0 0 0 0]$ and for  $X_3$  are $[0 1 1 1 1 0]$, $[1 0 1 1 0 1]$, $[1 1 0 1 1 1]$, $[1 1 1 1 0 0]$, and $[1 1 0 0 0 0]$. For $C$ they are $[0 1 1 1 1 0]$, $[1 0 1 1 0 0]$, $[1 1 0 1 1 1]$, $[1 1 1 1 0 1]$, and $1 1 0 0 0 0]$.
In each case we were able to find an $M^*(K_5\backslash e)$-minor (Appendix, Table A2).  

Lastly, from Table 1a we see that $X_1$ and $X_3$ have two non-isomorphic single-element extensions $Y_1$ and $Y_2$. Suppose $M$ is a coextension of $Y_1$ or $Y_2$. Then $M$ has rank 5 and 13 elements. If we add Type I rows, then $M\backslash  13$ is 3-connected, and if we add Type II or III rows, then $M\backslash 12$ is 3-connected, except when the row added is $[0 0 0 0 0 0 1 1]$. So only two matroids must be specifically checked for an $M^*(K_5\backslash e)$-minor. They are $Y_1$ with row $[0 0 0 0 0 0 1 1]$ and $Y_2$ with row $[0 0 0 0 0 0 1 1]$.   In both cases case the resulting matroid has an $M^*(K_5\backslash e)$-minor. Thus we may conclude that if $M$ has an $E_5$ or $D_2$-minor, then the rank of $M$ is at most 5. 
$\qed$
\bigskip


\section {Proof of Theorem 1.4} 

In this section we prove Theorem 1.4. We also give characterizations of the class of binary matroids with no prism-dual minor and the class of binary matroids with no prism and prims-dual minor. 

\bigskip
\noindent{\bf Proof of Theorem 1.4.} Theorem 1.3 establishes that the exceptional matroids in the class (i.e. matroids without an exact 3-seperation induced by $P_9$) are either the infinite families or $(P_\triangle (F_7, F_7)\backslash z)$ or have rank at most 5. Clearly $F_7$ and $F_7^*$ have no prism minor. Since $M^*(K_5\backslash e)$ has rank 5, all extensions of $F_7^*$ up to $PG(3, 2)$ are in the excluded minor class. These are shown in Table 1a.  

To complete the proof we must show that $R_{17}$ is the extremal rank 5 binary matroid with no $M^*(K_5\backslash e)$-minor. To do this we will show that if  $M$ is a rank-$5$ binary $3$-connected non-regular matroid with no $M^*(K_5\backslash e)$ and an $E_5$ or $D_2$-minor, then $M\cong R_{17}$ or its 3-connected restrictions (except $P_9^*$ because it has only 9 elements). 

Table 2 implies that the only columns that can be added to $E_5$ to obtain a matroid with no $M^*(K_5\backslash e)$-minor are those that give $A$, $B$, $C$. That is, columns $[0 0 1 0 1]$, $[0 0 1 1 0]$, $[0 1 0 1 1]$, $[0 1 1 0 0]$ $[1 0 0 1 1]$, $[1 1 0 0 1]$, $[1 1 1 0 1]$. It is straigntforward to check that adding all of these columns gives us a matroid isomorphic to $R_{17}$. The details are in Table 5.

From Claim 6 of Theorem 1.3, we see that besides $A$, $B$, and $C$, the matroid $Z$ is the only coextension of $D_2$ with no $M^*(K_5\backslash e)$-minor. As noted earlier, $Z$ is an extension of  $R_{10}$ and is obtained by adding any one of the columns $[0 1 0 1 1]$, $[0 1 1 0 1]$ $[1 0 1 0 1]$ $[1 0 1 1 0]$ $[1 1 0 1 0]$ or $[1 1 1 1 1]$ to $R_{10}$. The only other extension of $R_{10}$ is $A$. Adding all of the above six columns to $Z$ gives us $R_{17}\backslash e$. Adding one additional column (corresponding to extension $A$) gives us $R_{17}$.
 
One final matter must be checked. It may be possible for $R_{17}$ or one of its deletion-minors to be an extension of the graph $(K_5\backslash e)+edge$ or the cograph $(K'_{3,3})^*$. We must rule out this possibility. To do so, first observe from Table 3a that $M^*(K_5\backslash e)$ has two non-regular extensions $(K_5\backslash e)^*+edge$ and $(K'_{3,3})^*$. Second, observe that $E_5$ has no minor isomorphic to the prism graph or its dual. Third, Table 2 lists all the 3-connected deletion-minors of $A$, $B$, $C$, making it clear that they have no $M^*(K_5\backslash e)$-minor. Lastly, Table 5 gives the single-element extensions of $A$, $B$ and $C$ using columns $[0 0 1 0 1]$, $[0 0 1 1 0]$, $[0 1 0 1 1]$, $[0 1 1 0 0]$ $[1 0 0 1 1]$ $[1 1 0 0 1]$ and $[1 1 1 0 1]$ (the other columns give an $E_4$ minor, which has a prism minor). These columns give four 12-element extensions, $D$, $E$, $F$, and $G$. These 12-element matroids have no graphic nor cographic single-element deletions. So, they belong to the excluded minor class. Moreover, all their extensions 
will likewise have no graphic nor cographic single-element deletions. We are justified in addding all these columns to $E_5$ to get $R_{17}$. Hence proved.
$\qed$
\bigskip 

\tiny
 \begin{center}
\begin{tabular}{|c|p{15em}|c|}
\hline
\bf{Matroid} & {\bf Extension column} & {\bf Name}  \\  \hline \hline
$A$  & $\bf [0 0 1 1 0]$ $[0 1 1 0 0]$ $[1 0 0 1 1]$            & $D$    \\  \hline
& $\bf [0 1 0 1 1]$      & $E$ \\  \hline
& $\bf [1 1 0 0 1]$      & $F$ \\  \hline
& $\bf [1 1 1 0 1]$      & $G$ \\  \hline\hline

$B$ &$[0 0 1 0 1]$ $[0 0 1 1 0]$ $[0 1 0 1 1]$  $[0 1 1 0 1]$  & $D$ \\ \hline
& $[1 1 0 0 1]$  $[1 1 1 0 1]$    & $E$ \\  \hline\hline

$C$ & $[0 0 1 0 1]$ $[0 1 0 1 1]$ $[1 0 0 1 1]$ $[1 1 1 0 1]$  & $F$   \\  \hline
& $[0 0 1 1 0]$  $[0 1 1 0 0]$  & $G$ \\  \hline\hline

$D$  & $\bf [0 1 0 1 1]$ $[0 1 1 0 0]$ $[1 0 0 1 1]$ & $H$ \\ \hline
& $\bf[1 1 0 0 1]$ $[1 1 1 0 1]$     & $E$ \\ \hline\hline

$E$  & $[0 0 1 1 0]$ $[0 1 1 0 0]$ $[1 0 0 1 1]$ & $H$ \\ \hline
& $\bf [1 1 0 0 1]$   & $J$ \\ \hline
& $\bf [1 1 1 0 1]$   & $K$ \\ \hline\hline

$F$  & $[0 0 1 1 0]$ $[0 1 1 0 0]$ $[1 0 0 1 1]$ $[1 1 1 0 1]$ & $I$ \\ \hline
& $[0 1 0 1 1]$   & $J$ \\ \hline\hline

$G$  & $[0 0 1 1 0]$ $[0 1 1 0 0]$ $[1 0 0 1 1]$ $[1 1 0 0 1]$ & $I$ \\ \hline
& $[0 1 0 1 1]$   & $K$ \\ \hline\hline

$H$  & $\bf [0 1 1 0 0]$ $[1 0 0 1 1]$  & $L$ \\ \hline
& $\bf [1 1 0 0 1]$ $[1 1 1 0 1]$  & $M$ \\ \hline\hline

$I$  & $[0 1 0 1 1]$ $[0 1 1 0 0]$ $[1 0 0 1 1]$ $[1 1 1 0 1]$  & $M$ \\ \hline\hline

$J$  & $[0 0 1 1 0]$ $[0 1 1 0 0]$ $[1 0 0 1 1]$ $[1 1 1 0 1]$  & $M$ \\ \hline\hline

$K$  & $[0 0 1 1 0]$ $[0 1 1 0 0]$ $[1 0 0 1 1]$ $[1 1 0 0 1]$  & $M$ \\ \hline\hline

$L$  & $\bf [1 0 0 1 1]$   & $O$ \\ \hline
  & $\bf [1 1 0 0 1]$ $[1 1 1 0 1]$   & $P$ \\ \hline\hline

$M$  & $[0 1 1 0 0]$ $[1 0 0 1 1]$ $[1 1 1 0 1]$  & $P$ \\ \hline \hline

$O$  & $\bf [1 1 0 0 1]$ $[1 1 1 0 1]$   & $Q$ \\ \hline \hline

$P$  & $[1 0 0 1 1]$ $[1 1 1 0 1]$   & $Q$ \\ \hline \hline

$Q$  & $\bf [1 1 1 0 1]$   & $R$ \\ \hline \hline

\end{tabular}
\end{center}
\normalsize
 \begin{center} Table 5: All extensions of $E_5$ up to $R_{17}$ \end{center} 
\bigskip

Using Table 5, we can identify the internally 4-connected restrictions of $R_{17}$ as all, except $B$, $G$ and $K$. Among restrictions of $PG(3, 2)$ all except $K_5\backslash e$, $S_8$, $AG(3,2)$, $P_9$, $Z_4$, $D_1$ and $X_2$ are internally 4-connected. The next corollary follows immediately.

\bigskip
\noindent{\bf Corollary 3.1.} {\it$M$ is an internally $4$-connected binary matroid with no $M^*(K_5\backslash e)$-minor if and only if $M$ is isomorphic to $F_7$ or an internally $4$-connected restriction of $PG(3, 2)$ or $R_{17}$.} $\qed$
\bigskip

The above corollary is the main theorem in [5] by Mayhew and Royle.  The matroid they call $AG(3, 2)\times U_{1,1}$ is $R_{17}$. The five matroids they refer to are $B$, $G$, $K$, $D_1$, $X_2$. 
\bigskip

\noindent {Theorem 3.2.} {\it  Let $M$ be a binary matroid with no prism-minor.
\begin {enumerate}
\item [(i)] If $M$ is internally $4$-connected, then $M$ has rank at most $5$, and is isomorphic to a minor of $R_{17}$

\item [(i)] If $M$ is $3$-connected but not internally $4$-connected, and $M$ has an internally $4$-connected minor with at least $6$ elements that is not
isomorphic to $M(K_4)$, $F_7$, $F_7^*$ or $M(K_{3,3})$, then $M$ is isomorphic
to one of five matroids. \end {enumerate}}

\bigskip

Using Theorem 1.4 we can also identify members in the dual class (the class of binary matroids with no $M(K_5\backslash e)$-minor). However, we can only conclude that either $P_9$ or $P_9^*$ are 3-decomposers, instead of the stronger statment that ``$P_9$ is a 3-decomposer." This is because $E_7^*$ has no $M(K_5\backslash e)$-minor nor $P_9$-minor, but it does admit the 3-separation of its minor $P_9^*$. 

\bigskip
\noindent{\bf Theorem 3.3.} {\it Suppose $M$ is a $3$-connected binary matroid with no $M(K_5 \backslash e)$-minor. Then either $P_9$ or $P_9^*$ is a $3$-decomposer or $M$ is isomorphic to one of the following matroids: 
\begin {enumerate}
\item [(i)] $M^*(K_5\backslash e)$, $M(K_{3,3})$, $M^*(K_{3,3})$, or $M(W_r)$ for some $r\ge 3$; 
\item [(ii)]$Z_r$, $Z_r^*$, $Z_r\backslash b_r$, or $Z_r\backslash c_r$, for some $r\ge 4$; or 
\item [(iii)] $F_7$, $F_7^*$, $P_9$, $P_9^*$ 
\item [(iii)] $R_{17}^*$ or one of its contraction-minors.   $\qed$
\end{enumerate} }
\bigskip

\noindent {\bf Proof.} Suppose $M$ is a $3$-connected binary non-regular matroid  with no $M(K_5\backslash e)$-minor. If $M$ is regular, then $M$ is isomorphic to $M^*(K_5 \backslash e)$, $M(K_{3,3})$, $M^* (K_{3,3})$, $M(W_r)$ for some $r \ge 3$, or $R_{10}$. Therefore, suppose $M$ is non-regular. If $M$ has no $P_9$ or $P_9^*$-minor, then by Theorem 2.2 (Claim 1), $M$ is isomorphic to $F_7$, $F_7^*$, $Z_r$, $Z_r^*$, $Z_r\backslash b_r$, or $Z_r\backslash c_r$, for some $r\ge 4$. Thus we may assume $M$ has a $P_9$ or $P_9^*$-minor. Note that among the extensions of $P_9$, $D_2$ has an $M(K_5\backslash e)$-minor and $P_9$ is a 3-decomposer for rank-4 binary 3-connected matroids without a $D_2$-minor. So no further rank 4 matroids are int he class. It follows from Theorem 1.3(iv)   that $M$ is isomorphic to $R_{17}^*$ or its contraction minors. $\qed$
\bigskip

The next result follows from Theorems 1.3 and 3.3 and the fact that all the rank 5 exceptions and their duals have a minor isomorphic to $M(K_5\backslash e)$ or $M^*(K_5\backslash e)$ except $P_9$, $P_9^*$, and $E_5$.
\bigskip

\noindent{\bf Theorem 3.4.} {\it Suppose $M$ is a $3$-connected binary matroid with no $M(K_5\backslash e)$ or $M^*(K_5\backslash e)$-minor, then either $P_9$ is a $3$-decomposer for $M$ or  $M$ is isomorphic to $M(W_r)$ for some $r\ge 3$, $M(K_{3,3})$, $M^*(K_{3,3})$,  $Z_r$, $Z_r^*$, $Z_r\backslash b_r$, $Z_r\backslash c_r$, for some $r\ge 4$, $F_7$, $F_7^*$, $P_9$, $(P_\triangle (F_7, F_7)\backslash z)$, $(P_\triangle (F_7, F_7)\backslash z)^*$, $P_9^*$, $R_{10}$, or $E_5$.  }  $\qed$

\medskip

\small
\noindent {\bf Acknowledgement:} A year and half ago, the first author submitted a research proposal outlining matroid representability problems. The unknown reviewer indicated we should tackle this problem. The authors thank the reviewer for highlighting it. 
\normalsize
\bigskip

\noindent {\bf References}

\begin{enumerate}

\item  G. A. Dirac, Some results concerning the structure of graphs, {\it Canad. Math. Bull.} {\bf 6} (1963) 183-210. 

\item S. R. Kingan, Binary matroids without prisms, prism duals, and cubes, {\it Discrete Mathematics}, {\bf 152} (1996), 211-224.

\item S. R. Kingan, M. Lemos, Strong Splitter Theorem (submitted).

\item D. Mayhew, G. Royle, and G. Whittle, The internally 4-connected binary matroids with no $M(K_{3,3})$-minor. {\it Memoirs of the American Mathematical Society} {\bf 981}, American Mathematical Society, Providence, Rhode Island.

\item D. Mayhew, G. Royle, The internally 4-connected binary matroids with no $M(K_5\backslash e)$-minor (to appear) {\it Siam Journal on Discrete Mathematics.}

\item J. G. Oxley, The binary matroids with no 4-wheel minor, {\it Trans. Amer. Math. Soc.}  {\bf 154}  (1987), 63-75.

\item J. G. Oxley, {\it Matroid Theory}, (1992), Oxford University Press, New York.

\item N. Robertson and P. D. Seymour, Generalizing Kuratowski's Theorem, {\it Congressus Numerantium} {\bf 45} (1984) 129-138.

\item P. D. Seymour, Decomposition of regular matroids,  {\it J. Combin. Theory   Ser. B}  {\bf 28 }, (1980) 305-359.

\end {enumerate}

\vfill
\eject

\section {Appendix}

Table A1 lists the single-element coextensions of $A$, $B$, $C$, and $Z$.   Table A2 lists the single-element coextensions of $X_1$ and $X_3$.   
Type II and III are marked in red.  
\bigskip

\tiny
 \begin{center}
\begin{tabular}{|c|c|p{35em}|}
\hline
\bf{Matroid} & Name  & {\bf Coextension Row}  \\  \hline \hline
$A$ & coext 1 & \color {red} $[0 0 0 0 1 1]$  $[0 0 0 1 0 1]$ \color {black} $[0 0 1 0 1 0]$ $[0 1 1 0 1 0]$ $[1 0 1 1 1 1]$ $[1 1 1 0 0 1]$               \\  \hline

& coext 2 & $[0 0 0 1 1 0]$  $[1 1 0 0 1 1]$  $[1 1 0 1 0 1]$    \\  \hline

& coext 3 & $[0 0 0 1 1 1]$ \color {black} $[1 0 1 0 1 1]$ $[1 1 1 0 1 1]$ \\  \hline

& coext 4 & \color {black}$[0 0 1 0 0 1]$\color {black} $[0 1 0 1 1 0]$ \color {red}$[0 1 1 1 1 1]$ \color {black}      \\ \hline

& coext 5 & $[0 0 1 0 1 1]$ $[0 1 1 0 1 1]$ \color {black} $[1 0 0 1 1 1]$      \\ \hline

& coext 6 & $[0 0 1 1 0 0]$ $[0 1 1 1 0 0]$ \color {red}$[1 1 0 0 0 0]$      \color {black} \\ \hline

& coext 7 & $[0 0 1 1 0 1]$ \color {black} $[0 1 0 0 1 0]$ $[0 1 0 1 0 0]$  \color {black}  $[0 1 1 1 0 1]$  $[1 0 1 1 1 0]$  $[1 1 1 0 0 0]$     \\ \hline

& coext 8 & \color {black} $[0 0 1 1 1 0]$ $[0 1 1 0 0 0]$ \color {red} $[1 0 1 1 0 1]$ \color {black} $[1 1 0 0 1 0]$  $[1 1 0 1 0 0]$ \color {red} $[1 1 1 1 0 1]$      \\ \hline

& coext 9 & \color {black}  $[0 0 1 1 1 1]$ $[0 1 1 0 0 1]$ \color {black} $[1 0 0 0 1 1]$  $[1 0 0 1 0 1]$ \color {black}  $[1 0 1 0 1 0]$  $[1 1 1 0 1 0]$     \\ \hline

& coext 10 & \color {red} $[0 1 0 0 0 1]$ \color {black}  $[1 0 0 0 1 0]$ \color {black} $[1 0 0 1 0 0]$     \\ \hline

& coext 11 & \color {black}  $[0 1 0 0 1 1]$ $[0 1 0 1 0 1]$ $[1 0 0 1 1 0]$    \\ \hline

& coext 12 & $[0 1 0 1 1 1]$ \\ \hline

& coext 13 & \color {red} $[1 0 0 0 0 1]$ \color {black}  $[1 0 1 0 0 0]$ $[1 1 1 1 1 0]$    \\ \hline

& coext 14 & $[1 0 1 0 0 1]$ \color {red} $[1 1 0 1 1 0]$ \color {black} $[1 1 1 1 1 1]$     \\ \hline

 \hline

$B$ & coext 1  & \color {red}$[0 0 0 0 1 1]$ $[0 0 0 1 0 1]$ \color {black}  $[0 0 0 1 1 0]$  \color {red} $[0 0 1 0 0 1]$ \color {black}  $[0 0 1 0 1 0]$ \color {black} $[0 0 1 1 1 1]$ $[0 1 0 0 1 0]$ $[0 1 0 1 0 0]$ \color {black}  $[0 1 0 1 1 1 ]$ \color {black}  $[0 1 1 0 0 0]$ \color {black}  $[0 1 1 0 1 1]$ \color {red} $[0 1 1 1 1 0]$ \color {black}             \\  \hline

& coext 2 &  $[0 0 0 1 1 1]$ $[0 0 1 0 1 1]$ $[0 1 0 1 1 0]$ $[0 1 1 0 1 0]$            \\  \hline

& coext 3 &  $[0 0 1 1 0 0]$ \color {red} $[0 1 0 0 0 1]$ \color {black}  $[0 1 1 1 0 1]$           \\  \hline

& coext 4 &  $[0 0 1 1 0 1]$ \color {black}  $[0 0 1 1 1 0]$ $[0 1 0 0 1 1]$ $[0 1 0 1 0 1]$ $[0 1 1 0 0 1]$ \color {black}  $[0 1 1 1 0 0]$             \\   \hline

& coext 5 &  \color {red} $[1 0 0 0 0 1]$ \color {black} $[1 0 0 0 1 0]$ $[1 0 0 1 0 0]$ $[1 0 1 0 0 0]$ \color {red} $[1 0 1 1 0 1]$ \color {black} $[1 0 1 1 1 0]$ \color {red} $[1 1 0 0 0 0]$ \color {black} $[1 1 0 0 1 1]$ $[1 1 0 1 0 1]$ $[1 1 1 0 0 1]$ \color {red} $[1 1 1 1 0 0]$ \color {black} $[1 1 1 1 1 1]$             \\   \hline

& coext 6 & $[1 0 0 0 1 1]$ $[1 0 0 1 0 1]$ \color {black}  $[1 0 1 0 1 0]$ \color {black}  $[1 0 1 1 1 1]$ $[1 1 1 0 0 0]$ \color {black} $[1 1 1 0 1 1]$              \\   \hline

& coext 7 & \color {black}  $[1 0 0 1 1 0]$ \color {black} $[1 0 1 0 0 1]$ $[1 1 0 0 1 0]$ $[1 1 0 1 0 0]$ \color {red} $[1 1 0 1 1 1]$ \color {black} $[1 1 1 1 1 0]$            \\   \hline

& coext 8 & \color {black}  $[1 0 0 1 1 1]$ $[1 0 1 0 1 1]$ $[1 1 1 0 1 0]$             \\   \hline\hline

$C$ & coext 1 & \color {red} $[0 0 0 0 1 1]$ $[0 0 0 1 0 1]$ $[0 0 1 0 0 1]$  \color {black} $[0 0 1 1 1 1]$         $[0 1 0 0 1 0]$ $[0 1 0 1 0 0]$  $[0 1 1 0 0 0]$  \color {red} $[0 1 1 1 1 0]$ \color {black} $[1 0 0 0 1 0]$ $[1 0 0 1 0 0]$ $[1 0 1 0 0 0]$ $[1 0 1 1 1 0]$ $[1 1 0 0 1 1]$ $[1 1 0 1 0 1]$ $[1 1 1 0 0 1]$ $[1 1 1 1 1 1]$                   \\  \hline

& coext 2 & $[0 0 0 1 1 0]$  $[0 1 0 1 1 1]$     \\  \hline

& coext 3 & $[0 0 0 1 1 1]$  $[0 1 0 1 1 0]$  \color {black} $[1 0 0 1 1 0]$  \color {red} $[1 1 0 1 1 1]$ \color {black}  \\  \hline

& coext 4 & $[0 0 1 0 1 0]$  $[0 1 1 0 1 1]$    \\  \hline

& coext 5 & $[0 0 1 0 1 1]$  $[0 1 1 0 1 0]$ \color {black} $[1 0 1 0 1 0]$ $[1 1 1 0 1 1]$     \\  \hline

& coext 6 & $[0 0 1 1 0 0]$  $[0 1 1 1 0 1]$     \\  \hline

& coext 7 & $[0 0 1 1 0 1]$  $[0 1 1 1 0 0]$ \color {red} $[1 0 1 1 0 0]$ $[1 1 1 1 0 1]$ \color {black}  \\  \hline

& coext 8 & \color {black} $[0 0 1 1 1 0]$  $[0 1 0 0 1 1]$ $[0 1 0 1 0 1]$ $[0 1 1 0 0 1]$     \\  \hline

& coext 9 & \color {red} $[0 1 0 0 0 1]$ \color {black}   \\  \hline

& coext 10 & \color {red} $[1 0 0 0 0 1]$  $[1 1 0 0 0 0]$ \color {black}      \\  \hline

& coext 11 & $[1 0 0 0 1 1]$  $[1 0 0 1 0 1]$ $[1 0 1 1 1 1]$ $[1 1 1 0 0 0]$    \\  \hline

& coext 12 & \color {black}  $[1 0 0 1 1 1]$     \\  \hline

& coext 13 & $[1 0 1 0 0 1]$  $[1 1 0 0 1 0]$ $[1 1 0 1 0 0]$ $[1 1 1 1 1 0]$   \\  \hline

& coext 14 & \color {black}  $[1 0 1 0 1 1]$  $[1 1 1 0 1 0]$  \\  \hline\hline

$Z$ & coext 1 & \color {red} 
$[0 0 0 0 1 1]$ 
$[0 0 0 1 0 1]$ 
$[ 0 0 1 0 0 1 ]$
$[ 0 0 1 1 1 0 ]$
$[ 0 1 0 0 0 1 ]$
$[ 0 1 1 1 0 0 ]$
$[ 1 0 0 0 0 1 ]$
$[ 1 0 0 1 1 1 ]$
$[ 1 1 0 0 1 0 ]$
$[ 1 1 1 0 0 1 ]$ \color {black}
\\  \hline

& coext 2 & 
$[ 0 0 0 1 1 0 ]$
$[ 0 0 1 0 1 1 ]$
$[ 0 0 1 1 0 1 ]$
$[ 0 1 0 0 1 0 ]$
$[ 0 1 0 1 0 1 ]$
$[ 0 1 1 0 0 0 ]$
$[ 0 1 1 1 1 1 ]$
$[ 1 0 0 0 1 0 ]$
$[ 1 0 0 1 0 0 ]$
$[ 1 0 1 0 0 0 ]$
$[ 1 0 1 1 1 0 ]$
$[ 1 1 0 0 0 0 ]$
$[ 1 1 0 1 1 1 ]$
$[ 1 1 1 0 1 1 ]$
$[ 1 1 1 1 0 0 ]$
     \\  \hline

& coext 3 & 
$[ 0 0 0 1 1 1 ]$
$[ 0 0 1 0 1 0 ]$
$[ 0 0 1 1 0 0 ]$
$[ 0 1 0 0 1 1 ]$
$[ 0 1 0 1 0 0 ]$
$[ 0 1 1 0 0 1 ]$
$[ 0 1 1 1 1 0 ]$
$[ 1 0 0 0 1 1 ]$
$[ 1 0 0 1 0 1 ]$
$[ 1 0 1 0 0 1 ]$
$[ 1 0 1 1 1 1 ]$
$[ 1 1 0 0 0 1 ]$
$[ 1 1 0 1 1 0 ]$
$[ 1 1 1 0 1 0 ]$
$[ 1 1 1 1 0 1 ]$\\  \hline

& coext 4 & 
$[ 0 1 0 1 1 0 ]$
$[ 0 1 0 1 1 1 ]$
$[ 0 1 1 0 1 0 ]$
$[ 0 1 1 0 1 1 ]$
$[ 1 0 1 0 1 0 ]$
$[ 1 0 1 0 1 1 ]$
$[ 1 0 1 1 0 0 ]$
$[ 1 0 1 1 0 1 ]$
$[ 1 1 0 1 0 0 ]$
$[ 1 1 0 1 0 1 ]$
$[ 1 1 1 1 1 0 ]$
$[ 1 1 1 1 1 1 ]$    \\  \hline \hline

\end{tabular}
\end{center}
\normalsize
 \begin{center} Table A1: Single-element coextensions of $A$, $B$ and $C$ \end{center}

\bigskip

\tiny
 \begin{center}
\begin{tabular}{|c|c|p{35em}|}
\hline
\bf{Matroid} & {\bf  Name}  & {\bf Coextension Row}  \\  \hline \hline
$X_1$ & coext 1 &  \color {red}  $[0 0 0 0 0 1 1]$ $[0 0 0 0 1 0 1]$ \color {black}  $[0 0 0 0 1 1 0]$ \color {red}  $[0 0 0 1 0 0 1]$ \color {black}  $[0 0 0 1 0 1 0]$ $[0 0 0 1 1 0 0]$ $[0 0 1 0 0 1 1]$ $[0 0 1 1 1 0 0 ]$ $0 1 0 0 1 1 0]$ $[0 1 0 1 0 0 1]$ $[0 1 1 0 1 0 1]$ $[0 1 1 1 0 1 0]$
 \\  \hline

& coext 2 & $[0 0 0 1 1 1]$ $[0 0 1 0 1 1]$ $[0 0 0 1 1 0 1]$ $[0 0 0 1 1 1 0]$ $[1 0 0 1 1 1 1]$ $[1 0 1 0 0 1 1]$ $[1 1 0 0 1 1 0]$ $[1 1 1 0 1 0 1]$    \\  \hline

& coext 3 & $[0 0 0 1 1 1 1]$ \\  \hline

& coext 4 & \color {red}  $[0 0 1 0 0 0 1]$ \color {black}  $[0 0 1 0 0 1 0]$ $[0 0 1 0 1 0 0]$ $[0 0 1 1 0 0 0]$ \color {red}  $[0 1 0 0 0 0 1]$ \color {black}  $[0 1 0 0 0 1 0]$ $[0 1 0 0 1 0 0]$ $[0 1 0 1 0 0 0]$ $[0 1 1 0 1 1 1]$ $[0 1 1 1 0 1 1]$ $[0 1 1 1 1 0 1]$ \color {red}  $[0 1 1 1 1 1 0]$ \color {black}  $[1 0 0 0 1 0 1]$ $[1 0 0 1 0 0 1]$ $[1 0 0 1 1 0 0]$ $[1 0 1 0 0 0 0]$ $[1 0 1 0 1 1 0]$ $[1 0 1 1 0 1 0]$ $[1 1 0 0 0 0 0]$ $[1 1 0 0 0 1 1]$ $[1 1 0 1 0 1 0]$ $[1 1 1 1 0 0 1]$ $[1 1 1 1 1 0 0]$ $[1 1 1 1 1 1 1]$    \\ \hline

& coext 5 & $[0 0 1 0 1 0 1]$ $[0 0 1 0 1 1 0]$ $[0 0 1 1 0 0 1]$ $[0 0 1 1 0 1 0]$ $[0 1 0 0 0 1 1]$ $[0 1 0 0 1 0 1]$ $[0 1 0 1 0 1 0]$ $[0 1 0 1 1 0 0]$ $[0 1 1 0 0 1 1]$ $[0 1 1 0 1 1 0]$ $[0 1 1 1 1 0 0 1]$ $[0 1 1 1 1 0 0]$      \\ \hline

& coext 6 & $[0 0 1 0 1 1 1]$ $[0 0 1 1 0 1 1]$ $[0 0 1 1 1 0 1]$ $[0 0 1 1 1 1 0]$ $[0 1 0 0 1 1 1]$ $[0 1 0 1 0 1 1]$ $[0 1 0 1 1 0 1]$ $[0 1 0 1 1 1 0]$ $[0 1 1 0 0 0 1]$ $[0 1 1 0 0 1 0]$ $[0 1 1 0 1 0 0]$ $[0 1 1 1 0 0 0]$ $[1 0 0 0 0 1 1]$ $[1 0 0 0 1 1 0]$ $[1 0 0 1 0 1 0]$ $[1 0 1 0 1 0 1]$ $[1 0 1 1 0 0 1]$ $[1 0 1 1 0 0 1]$ $[1 0 1 1 1 1 1]$ $[1 1 0 0 1 0 1]$ $[1 1 0 1 1 0 0]$ $[1 1 0 1 1 1 1]$ $[1 1 1 0 0 0 0]$ $[1 1 1 0 0 1 1]$ $[1 1 1 0 1 1 0]$\\ \hline

& coext 7 & $[ 0 0 1 1 1 1 1]$ $[0 1 0 1 1 1 1]$ $[0 1 1 0 0 0 0]$     \\ \hline

& coext 8 & \color {red}  $[1 0 0 0 0 0 1]$ \color {black}  $[1 0 0 0 1 0 0]$ $[1 0 0 1 0 0 0]$ $[1 0 1 0 1 0 0 0]$ $1 0 1 1 0 0 0]$ $[1 0 1 1 1 1 0]$ $[1 1 0 0 0 1]$ \color {red}  $[1 1 0 1 0 0 0]$ \color {black}  $1 1 0 1 0 1 1]$   $[1 1 1 1 0 0 0]$ \color {red}  $[1 1 1 1 0 1 1]$ \color {black}  $[1 1 1 1 1 1 0]$    \\ \hline

& coext 9 & $[1 0 0 0 1 0]$ \color {red}  $[1 0  1 1 1 0 1]$ \color {black}  $[1 1 0 1 1 0 1]$ $[1 1 1 0 0 1 0]$ \\ \hline

& coext 10 & $[1 0 0 0 1 1 1]$ $[1 0 0 1 0 1 1]$ $[1 0 0 1 1 1 0]$ $[1 0 1 0 0 0 1]$ $[1 0 1 0 1 1 1]$ $[1 0 1 1 0 1 1]$ $[1 1 0 0 1 0 0]$ $[1 1 0 0 1 1 1]$ $[1 1 0 1 1 1 0]$ $[1 1 1 0 0 0 1]$ $[1 1 1 0 1 0 0]$ $[1 1 1 0 1 1 1]$    \\ \hline

& coext 11 & $[1 0 0 1 1 0 1]$ $[1 0 1 0 0 1 0]$ $[1 1 0 0 0 1 0]$ $[1 1 1 1 1 0 1]$    \\ \hline\hline

$X_3$ & coext 1  &    \color {red} $[0 0 0 0 0 1 1]$  \color {black} $[0 0 0 0 1 1 0]$ \color {red}  $[0 0 1 0 0 0 1]$ \color {black}  $[0 1 0 0 1 1 1]$ $[0 1 0 1 1 1 0]$ $[0 1 1 0 0 0 0]$ $[0 1 1 0 1 0 1]$ $[1 0 1 0 1 1 1]$ $[1 0 1 1 0 1 1]$ $[1 1 0 1 1 0 1]$ $[1 1 1 0 0 1 1]$        \\  \hline

& coext 2 &  \color {red} $[0 0 0 0 1 0 1]$ $[0 0 0 1 0 0 1]$ \color {black}  $[0 0 0 1 1 0 0]$ $[0 0 1 1 1 1 0]$ \color {red}  $[0 1 0 0 0 0 1]$ \color {black} $[0 1 0 1 0 0 0]$ $[0 1 1 1 0 1 0]$ \color {red} $[0 1 1 1 1 1 1]$ \color {black}  $[1 0 1 0 1 0 0]$ $[1 0 1 1 0 0 0]$ $[1 1 0 0 0 1 0 ]$ $[1 1 1 1 1 0 0]$          \\  \hline

& coext 3 &  $[0 0 0 0 1 1 1]$ $[0 0 0 1 0 1 1]$ $[0 0 0 1 1 1 0]$ $[0 1 0 1 1 1 1]$ $[0 1 1 0 0 0 1]$ $[1 0 1 0 0 1 1]$          \\  \hline

& coext 4 &  $[0 0 0 1 1 0 1]$ $[0 0 1 1 1 1 1]$ $[1 1 0 0 1 1 0]$ $[1 1 1 0 1 0 0]$ \\   \hline

& coext 5 &  $[0 0 0 1 1 1 1]$             \\   \hline

& coext 6 & $[0 0 1 0 0 1 0]$ $[1 1 0 1 0 1 1]$ $[1 1 1 1 0 0 1]$    \\   \hline

& coext 7 &$[0 0 1 0 0 1 1]$ $[1 0 0 0 1 1 1]$ $[1 0 0 1 0 1 1]$ $[1 0 0 1 1 1 0]$ $[1 1 0 1 1 1 1]$ $[1 1 1 0 0 0 1]$            \\   \hline

& coext 8 & $[0 0 1 0 1 0 0]$ $[0 0 1 1 0 0 0]$ $[0 1 0 0 0 1 0]$ $[0 1 1 1 1 0 0]$ $[1 0 0 0 1 0 1]$ $[1 0 0 1 0 0 1]$ $[1 0 0 1 1 0 0]$ $[1 0 1 1 1 1 0]$ $[1 1 0 0 0 0 1]$ \color {red}  $[1 1 0 1 0 0 0]$ $[1 1 1 1 0 1 0]$ \color {black}  $[1 1 1 1 1 1 1]$      \\   \hline

& coext 9 & $[0 0 1 0 1 0 1]$ $[0 0 1 1 0 0 1]$ $[0 1 0 0 0 1 1]$ $[0 1 0 1 0 1 0]$ $[0 1 1 1 0 0 0]$ $[0 1 1 1 1 0 1]$ $[1 0 1 0 1 1 0]$ $[1 0 1 1 0 1 0]$ $[1 1 0 0 1 0 1]$ $[1 1 0 1 1 0 0]$ $[1 1 1 0 0 1 0]$ $[1 1 1 0 1 1 1]$   \\   \hline

& coext 10 & $[0 0 1 0 1 1 0]$ $[0 0 1 1 0 1 0]$ $[0 1 0 0 1 0 1]$ $[0 1 0 1 1 0 0]$ $[0 1 1 0 0 1 0]$ $[0 1 1 0 1 1 1]$ $[1 0 1 0 1 0 1]$ $[1 0 1 1 0 0 1]$ $[1 1 0 0 0 1 1]$ $[1 1 0 1 0 1 0]$ $[1 1 1 1 0 0 0]$ $[1 1 1 1 1 0 1]$     \\   \hline

& coext 11 & $[0 0 1 0 1 1 1]$ $[0 0 1 1 0 1 1]$  $[0 1 1 0 0 1 1]$ $[1 0 0 0 0 1 1]$ $[1 0 0 0 1 1 0]$ $[1 0 0 1 0 1 0]$ $1 0 1 0 0 0 1]$ $1 1 0 0 1 1 1]$ $[1 1 0 1 1 1 0$ $[1 1 1 0 0 0 0]$ $[1 1 1 0 1 0 1]$     \\   \hline

& coext 12 & $[0 0 1 1 1 0 0]$ \color {red} $[1 0 0 0 0 0 1]$ \color {black} $[1 0 0 0 1 0 0]$ $[1 0 0 1 0 0 0]$ $1 1 0 0 0 0 0]$ $[1 1 1 1 1 1 0]$     \\   \hline

& coext 13 &  $[0 0 1 1 1 0 1]$ $[1 1 0 0 1 0 0]$ $[1 1 1 0 1 1 0]$    \\   \hline

& coext 14 & $[0 1 0 0 1 0 0]$ $[0 1 1 0 1 1 0]$ \color {red}  $[1 0 1 1 1 0 1]$ \color {black}     \\   \hline

& coext 15 & $[0 1 0 0 1 1 0]$ $[0 1 1 0 1 0 0]$ $[1 0 0 1 1 0 1]$ $[1 0 1 1 1 1 1]$    \\   \hline

& coext 16 & $[0 1 0 1 0 0 1]$ $[0 1 1 1 0 1 1]$ $[1 0 0 0 0 1 0 ]$ $[1 0 1 0 0 0 0]$    \\   \hline

& coext 17 & $[0 1 0 1 0 1 1]$ $[0 1 1 1 0 0 1]$ $[1 0 1 0 0 1 0]$     \\   \hline

& coext 18 & $[1 0 0 1 1 1 1]$     \\   \hline\hline

\end{tabular}
\end{center}
\normalsize
 \begin{center} Table A2: Single-element coextensions of $X_1$ and $X_3$ \end{center}

\end {document}